\documentclass[a4paper]{article}
\usepackage[utf8]{inputenc}
\usepackage{fullpage}
\usepackage[centertags]{amsmath}
\usepackage{color}
\usepackage{amssymb}
\usepackage{amsthm}
\usepackage{dsfont}
\usepackage[all]{xy}
\usepackage{mathrsfs}
\usepackage{upgreek}
\usepackage{tikz-cd}
\usepackage[colorlinks=true, linkcolor=blue, urlcolor=red, citecolor=magenta]{hyperref}
\usepackage{mathtools}
\usepackage{microtype}
\usepackage{lmodern}
\usepackage[nottoc]{tocbibind}
\usepackage{extarrows}
\usepackage{enumitem}
\usepackage[capitalise]{cleveref}
\crefname{enumi}{}{}
\usepackage{stmaryrd}
\usepackage{thmtools}

\newtheorem{theorem}{Theorem}[section]
\numberwithin{equation}{subsection}
\newtheorem{lemma}[theorem]{Lemma}
\newtheorem{proposition}[theorem]{Proposition}
\newtheorem{corollary}[theorem]{Corollary}

\theoremstyle{definition} \newtheorem{definition}[theorem]{Definition}
\newtheorem{example}[theorem]{Example}
\newtheorem{conjecture}[theorem]{Conjecture}
\newtheorem{assumption}[theorem]{Assumption}

\theoremstyle{remark}
\newtheorem{remark}[theorem]{Remark}

\newcommand{\nd}{\noindent}

\newcommand{\dR}{{\mathds R}}
\newcommand{\dC}{{\mathds C}}
\newcommand{\dQ}{{\mathds Q}}

\newcommand{\dZ}{{\mathds Z}}

\newcommand{\dL}{{\mathbb L}}
\newcommand{\dA}{{\mathbb A}}

\newcommand{\bD}{{\mathbb D}}

\newcommand{\cD}{\mathscr{D}}
\newcommand{\cE}{\mathcal{E}}
\newcommand{\cF}{\mathcal{F}}

\newcommand{\cM}{\mathcal{M}}

\newcommand{\cO}{\mathcal{O}}

\newcommand{\cR}{\mathcal{R}}

\newcommand{\cT}{\mathcal{T}}

\newcommand{\ft}{\mathfrak{t}}

\newcommand{\T}{\textstyle}

\DeclareMathOperator{\SL}{\textup{SL}}

\DeclareMathOperator{\pr}{\textup{pr}}

\DeclareMathOperator{\gr}{\textup{gr}}

\newcommand{\Kl}{\mathrm{Kl}}

\newcommand{\MHM}{\textup{MHM}}
\newcommand{\MTM}{\textup{MTM}}

\newcommand{\EMHM}{\textup{EMHM}}

\newcommand{\mbg}{\mathds{G}}

\newcommand{\mfg}{\mathfrak{g}}

\newcommand{\p}{\partial}

\DeclareMathOperator{\rk}{rk}

\DeclareMathOperator{\sheafHom}{\mathscr{H}\kern -3pt\textit{om}\kern 1pt}
\DeclareMathOperator{\sheafDer}{\mathscr{D}\kern -1pt\textit{er}\kern 1pt}

\newcommand{\C}{\mathds{C}}

\newcommand{\G}{\mathbb{G}}

\newcommand{\g}{\mfg}

\newcommand{\dif}{\mathrm{d}}

\newcommand{\al}{\alpha}

\newcommand{\bb}[1]{\mathbb{#1}}
\newcommand{\rr}[1]{\mathrm{#1}}
\newcommand{\cc}[1]{\mathcal{#1}}
\newcommand{\scr}[1]{\mathscr{#1}}

\makeatletter
\newcommand{\subjclass}[2][2020]{  \let\@oldtitle\@title  \gdef\@title{\@oldtitle\footnotetext{#1 \emph{Mathematics subject classification.} #2}}}
\newcommand{\keywords}[1]{  \let\@@oldtitle\@title  \gdef\@title{\@@oldtitle\footnotetext{\emph{Key words.} #1.}}}
\makeatother

\begin{document}
\title{\Large Irregular Hodge numbers of Frenkel--Gross connections}
\author{Yichen Qin, Christian Sevenheck, and Peter Spacek}


\maketitle

\begin{abstract}
  Frenkel and Gross constructed a family of connections on $\bb{P}^1\backslash\{0,\infty\}$, for almost simple groups $\check{G}$ and their representations. In this article, we calculate the irregular Hodge numbers of these Frenkel--Gross connections, and, as an application, we prove a conjecture of Katzarkov--Kontsevich--Pantev for mirror Landau-Ginzburg models of minuscule homogeneous spaces.
\end{abstract}

\tableofcontents

\section{Introduction} \label{sec:intro}
The study of irregular singularities of systems of linear differential equations has gained increasing importance across various mathematical fields, such as Hodge theory or arithmetics. Since the differential equations underlying variations of Hodge structures (VHS) can only have regular singularities, we need a broader framework to discuss the Hodge theoretic properties of irregular connections.

A famous example is the Bessel differential equation $(t\partial_t)^{n+1}-t$. The associated connection, called the \textit{Bessel connection} (or \emph{Kloosterman connection}), is the connection
\begin{equation}\label{eq:Bessel}
            \nabla=\dif+\begin{pmatrix}
         0 &  &  &   t \\
           1 & \ddots & &  \\
           & \ddots & \ddots &  \\
            &  & 1 & 0 \end{pmatrix}
        \frac{\dif t}{t}
\end{equation}
    on the trivial bundle $\cO^{n+1}_{\G_m}$, denoted by $\rr{Be}_{n+1}$. It has a regular singularity at $0$, and an irregular singularity at $\infty$ of slope $\frac{1}{n+1}$. When $n=1$, Deligne showed that $\rr{Be}_{2}$ cannot underlie a variation of Hodge structures \cite[\S8]{Deligne}, and hence, there does not exist any Hodge filtration on $\rr{Be}_{2}$. Instead, Deligne suggested and studied the concept of irregular Hodge filtrations on irregular connections \cite{Deligne}, which was further developed by Esnault, Kontsevich, Sabbah, and Yu in \cite{Sabbah-2010-laplace, Yu14, Sabbah-Yu-2015, Esnault-Sabbah-Yu}. For example, the irregular Hodge filtration on $\rr{Be}_{n+1}$ jumps at $0,1,\dots,n$, and the irregular Hodge numbers are all $1$, see e.g.~\cite{CDS21}.

Since then, a solid framework of the theory of irregular Hodge structures and irregular Hodge modules has been developed in \cite{Kontsevich-Soibelman, Sabbah-Irreg-book,FSY22, Mochizuki}. We will use the fundamental results of these papers in an essential way here.
Rather than using the term irregular Hodge modules, we adopt the terminology of rescalable integrable mixed twistor $\scr{D}$-modules, as introduced in  \cite{Mochizuki}, cf.~\cref{sec:EMHM}. In particular, the underlying $\cD$-module of any rescalable integrable mixed twistor $\scr{D}$-module comes equipped with an irregular Hodge filtration.

The purpose of this paper is to study irregular Hodge theory for certain rigid irregular connections on principal bundles. More precisely,
given any almost simple complex algebraic group $\check{G}$,
we are concerned with a family of connections (denoted $\nabla_{\check{G}}$) on the trivial $\check{G}$-bundle over $\bb{P}^1\backslash\{0,\infty\}$
that has been constructed by
Frenkel and Gross in \cite{FrenkelGross}, cf.~\eqref{eq:FG-conn}.  These connections are regular singular at $0$ and irregular singular at $\infty$. They are cohomologically rigid \cite[Thm.~1]{FrenkelGross} and expected to be physically rigid. The physical rigidity is not known to hold in general, but has been shown when $\check{G}$ is of adjoint type \cite[Thm.~1]{Yi22} or simply connected \cite[Thm\,1.1.1]{Jakob-Hohl}.
By considering specific representations of $\check{G}$, one can recover examples such as the Bessel connection and certain hypergeometric connections, cf.~\cref{exe:An}.

Our main result establishes the existence of a rescalable integrable mixed twistor $\scr{D}$-module associated to Frenkel–Gross connections and provides an explicit computation of the corresponding irregular Hodge numbers.

In the statement of the following theorem, we denote by $\rho=2\rho/2$ the half sum of positive coroots of $\check{G}$, which can also be considered as the half sum of positive roots of $G$. We recall a few more facts and notations around algebraic groups at the end of this introduction.

\begin{theorem}\label{thm:main}
    For an almost simple reductive group $\check G$ over $\dC$, the Frenkel--Gross connection $\nabla_{\check G}$ underlies a unique rescalable integrable mixed twistor $\scr{D}$-module $\nabla_{\check G}^{\rr{H}}$. The corresponding irregular Hodge numbers are determined\footnote{See \cref{def:FiltGivenCocarac}  for a precise definition.} by $\rho=2\rho/2$, the half sum of positive coroots of $\check{G}$. Concretely, for each representation $V$ of $\check{G}$, the irregular Hodge numbers $h_{\rr{irr}}^\alpha$ of $\nabla_{\check{G}}(V)$ are the dimensions of the eigenspaces for the eigenvalues $2\alpha$ of $V$ under the $\bb{G}_m$-action $ \bb{G}_m\xrightarrow{2\rho} \check G\to \rr{GL}(V)$.
\end{theorem}
 The irregular Hodge filtration (as opposed to the irregular Hodge numbers, i.e.~the dimension of the graded pieces of this filtration) is not a priori controlled by our main result. However, we formulate in \cref{sec:RMod} a conjecture that would lead to an alternative proof of our main result and would additionally imply that the irregular Hodge filtration itself is also determined by the cocharacter $\rho$ (in the sense of \cref{def:FiltGivenCocarac}).

\paragraph{Frobenius slopes.}
     One of Deligne’s motivations for studying irregular Hodge filtrations is to obtain lower bounds for $p$-adic slopes of exponential sums. Let us mention an arithmetic application of \cref{thm:main} for the case of Frenkel--Gross connections, as mentioned in \cite[Rem. 13.9]{Lam-Templier}. In \cite[Thm. 5.3.2]{XZ22}, Xu and Zhu equipped Frenkel--Gross connections with Frobenius structures and studied their $p$-adic slopes---i.e.~$p$-adic valuations of eigenvalues of Frobenius on Frenkel--Gross connections. The traces of the Frobenius actions are called \textit{Kloosterman sums (for reductive groups)}, which are sums of Frobenius eigenvalues. They showed that the $p$-adic slopes are determined by $\rho$ (except for finitely many points). In combination with this result, \cref{thm:main} implies that the $p$-adic slopes agree with the irregular Hodge numbers. In particular, when $\check{G}=\rr{SL}_{n+1}$, we recover a theorem of Sperber \cite{Sperber}.

\paragraph{Mirror symmetry.}
    Another motivation for studying irregular Hodge filtration arises from mirror symmetry. As mirrors of Fano varieties one considers Landau--Ginzburg models, namely, pairs $(Y,w)$ consisting of a quasi-projective complex variety and a regular function $w\colon Y\to \bb{A}^1$. A great deal of work has been carried out over the last decades to establish such mirror correspondences, see, e.g., \cite{Iritani} and \cite{RS10, RS12} for the toric case, and, on the other hand, \cite{Rietsch} and \cite{Lam-Templier} for the case of homogeneous spaces.

    In \cite{KKP2}, Katzarkov--Kontsevich--Pantev proposed three kinds of (conjecturally identical) Hodge numbers of $(Y,w)$, including the irregular Hodge numbers. (Notice, however, that one of these conjectured equalities is incorrectly stated, as can be seen by looking at the case of the mirror of projective spaces.) In particular, this yields a conjectured relationship comparing the Hodge numbers for Fano varieties and the irregular Hodge numbers for their Landau--Ginzburg models. As an application of \cref{thm:main}, we confirm this conjecture in \cref{thm:Kontsevich} for Fano varieties that are minuscule homogeneous spaces.

\paragraph{Idea of the proof of \texorpdfstring{\cref{thm:main}}{}.}

In order to orient the reader, let us sketch informally the strategy to show our main result \cref{thm:main}.
The first ingredient is the fact that the Frenkel--Gross connections ${\nabla}_{\check{G}}(V)$ are of (exponentially) geometrical origin. In fact, by a theorem of Zhu \cite{Zhu17}, they are eigen $\scr{D}$-modules of Hecke operators, constructed via the method of Heinloth--Ng\^o--Yun \cite{HeinlothNgoYun}. Using this strategy, we equip
in \cref{sec:EMHM} the connection ${\nabla}_{\check{G}}(V)$ with the structure of a rescalable integrable pure twistor $\scr{D}$-module of weight $0$. Furthermore, we show that there is a unique lifting of $\nabla_{\check{G}}$ to a tensor functor $\nabla_{\check{G}}^{\rr{H}}$ from the category of finite dimensional representations of $\check{G}$ to the category of rescalable integrable mixed twistor $\scr{D}$-modules on $\bb{G}_{m}$ (see the discussion before formula \eqref{eq:EMHM-tensor}).

Via the functoriality of $\nabla_{\check{G}}^{\rr{H}}$ (see \cref{sec:functoriality}), we reduce (Lemmas~\ref{lem:reduction1} and~\ref{lem:reduction2}) the calculation of Hodge numbers to that of $\nabla_{\check{G}}^{\rr{H}}(V)$ for simple groups $\check{G}$ and $V=\rr{Ad}(\check{G})$, the adjoint representation of $\check{G}$.

The next step is to express the Frenkel--Gross connections $\nabla_{\check{G}}(V)$ as the Fourier transforms of some regular singular connections. By the stationary phase principle (cf. \cite{sabbah2008explicit}), the Fourier transform of a regular $\scr{D}_{\bb{A}^1_\tau}$-module is a $\scr{D}_{\bb{A}^1_t}$-module (where $\dA^1_t$ is the dual affine line of $\dA^1_\tau$) which is smooth on $\bb{G}_{m,t}$, regular singular at $0$, and irregular singular at $\infty$ of slope $0$ or $1$. Frenkel and Gross showed that $\nabla_{\check{G}}(V)$ has slope $0$ or $1/h$ where $h$ is the Coxeter number of $\check{G}$. Therefore, we rather work with $\widetilde{\nabla}_{\check{G}}(V):=[h]^+\nabla_{\check{G}}(V)$, which has the same irregular Hodge numbers as $\nabla_{\check{G}}(V)$. Using the geometrical presentation of Frenkel--Gross connections mentioned above, we express $\widetilde{\nabla}_{\check{G}}(V)$ as the Fourier transform of some regular $\scr{D}_{\bb{A}^1_\tau}$-modules in \cref{prop:HNY-D-module}.

By a (filtered) stationary phase formula \cite{sabbah-hypergeometric}, to calculate the irregular Hodge numbers of $\widetilde{\nabla}_{\check{G}}(V)$, it suffices to calculate the Hodge numbers of the nearby cycle at infinity of $\rr{FT}^{-1}(j_+\nabla_{\check{G}}^{\rr{H}}(V))$, where $j:\G_{m,t}\hookrightarrow\dA^1_t$. In general, the inverse Fourier transform of Frenkel--Gross connections involves some intersection complexes, which makes the Hodge numbers tricky to handle. However, the nilpotent part of the monodromy of the Frenkel--Gross connections at $0$ is quite special, namely principally nilpotent. Through the stationary phase principle mentioned above, we can easily determine the monodromy of the inverse Fourier transform $\rr{FT}^{-1}\widetilde{\nabla}_{\check{G}}(V)$ at $\infty$ using representation-theoretic information of $V$. In particular, when the nilpotent part of the monodromy operator can be decomposed into Jordan blocks of different sizes, we calculate in \cref{thm:minuscule} the formula of Hodge numbers of the nearby cycle of $\rr{FT}^{-1}(j_+\nabla_{\check{G}}^{\rr{H}}(V))$ at $\infty$. Since this condition holds for adjoint representations, we conclude the proof of our main result.

\paragraph{Plan of the paper}
The content of the paper is organized as follows:
In \cref{sec:FGconn}, we collect basic properties of Frenkel--Gross connections. In \cref{sec:IrregHodge}, we recall preliminaries on rescalable integrable mixed twistor $\scr{D}$-modules and equip Frenkel--Gross connections with such structures. Moreover, we study the inverse Fourier transforms of $\widetilde{\nabla}_{\check{G}}$ in \cref{sec:inverse-fourier} and prove \cref{thm:main} in \cref{sec:proof}. In \cref{sec:KKP-conj}, we verify a conjecture of Katzarkov--Kontsevich--Pantev for minuscule flag varieties. In \cref{sec:exe}, we provide concrete calculations of irregular Hodge numbers for certain examples.
In \cref{sec:RMod}, we make a conjecture of the shape of the underlying $\scr{R}$-modules of $\nabla_{\check{G}}^{\rr{H}}$.

\medskip

We finish this introduction by fixing some notation that is used throughout this article:

\paragraph{$\cD$-modules.} Throughout the paper, we work with algebraic $\cD$-modules. One may consult \cite{Hotta} as a reference for the notation and results we use. For a smooth variety $X$ over $\dC$, we write $D^b_{qc}(\cD_X)$ resp. $D^b_{h}(\cD_X)$ for the derived category of quasi-coherent resp. of holonomic $\cD_X$-modules. For a morphism $f:X\rightarrow Y$ we have functors
$f_+ \colon D_{qc}^b(\mathscr D_X) \to D_{qc}^b(\mathscr D_Y)$ and $f^+ \colon D_{qc}^b(\mathscr D_Y) \to D_{qc}^b(\mathscr D_X)$. Moreover, there is the duality functor $\bD \colon D^b_h(\cD_X)\to
D^b_h(\cD_X)$ that respects $\textup{Mod}_h(\cD_X)$. We also define functors
$f_\dag:=\bD\circ f_+ \circ \bD$ and $f^\dag:=\bD\circ f^+ \circ \bD$.

We will use at various place the Fourier transformation functor for algebraic $\cD$-modules. The most general definition is when we are starting with a vector bundle $E\to X$, where $X$ is still smooth. Then we have functors $\mathrm{FT}^{\pm 1}_X:D^b(\cD_E)\to D^b(\cD_{E^\vee})$, where $E^\vee$ is the dual bundle. They are defined by
$\mathrm{FT}^{\pm 1}(N):=p_{2,+}(p_1^+ N \otimes \cE^{\pm\varphi})$, where $p_1:E\times E^\vee\rightarrow E$ and $p_2:E\times E^\vee\rightarrow E^\vee$ are the projections and where $\varphi\in\cO_{E\times E^\vee}$ is the canonical pairing. The functors $\mathrm{FT}^{\pm 1}_X$ are exact (i.e. they send $\textup{Mod}(\cD_E)$ to $\textup{Mod}(\cD_{E^\vee})$). Special cases of this definition that we will use are when $X=\{pt\}$ (in which case we write $\mathrm{FT}^{\pm 1}:=\mathrm{FT}^{\pm 1}_{pt}$) or when $E\cong \dA^1\times X$ is the trivial line bundle over $X$ (this case is called partial Fourier transformation with respect to the $\dA^1$-coordinate).

We consider the abelian category $\MHM(X)$ of algebraic
$\dQ$-mixed Hodge modules on $X$ (as defined in \cite{Saito1, SaitoMHM}) and we write $D^b \MHM(X)$ for its bounded derived category.
$\MHM(X)$ comes with a forgetful functor, sending a mixed Hodge module to its underlying regular holonomic $\scr{D}_X$-module. Furthermore, this forgetful functor induces a forgetful functor on the bounded derived categories $D^b(\mathrm{MHM}(X)) \to D^b_h(\scr{D}_X)$.
Let $f\colon X\rightarrow Y$ by any morphism of algebraic varieties, then the functors $f_+, f_\dag$ resp.\ $f^\dag[\dim(Y)-\dim(X)], f^+[\dim(X)-\dim(Y)]$ on $D^b_{h}(\mathscr{D}_X)$ resp.\
$D^b_{h}(\mathscr{D}_Y)$ lift to functors
$$
f_*, f_!\colon  D^b\MHM(X)\rightarrow D^b\MHM(Y)
\quad\quad
\textup{resp.}
\quad\quad
f^*, f^!\colon  D^b\MHM(Y)\rightarrow D^b\MHM(X).
$$
We also denote by $\bD$ the functor on $D^b \MHM(X)$ which lifts the above defined holonomic duality functor on $D^b_h(\cD_X)$.
We also consider $\dR$-mixed Hodge modules (as in \cite[\S 13.5]{Mochizuki-MTM}), and we denote the corresponding abelian category by $\MHM(X,\dR)$ and by
$D^b \MHM(X,\dR)$ the corresponding bounded derived category.

A major part of our considerations will use the category of integrable mixed twistor modules as constructed by Mochizuki (see \cite{Mochizuki-MTM}). We will postpone reminders on the properties of this category to \cref{sec:MTM}, where the motivation for using this category becomes clear.

\paragraph{Simple groups.} For an almost simple reductive group $G$ over $\C$, we fix a choice of a maximal torus $T$ and a Borel subgroup $B$ containing $T$. The \emph{root datum} of $G$ is the data $(X_*,\Phi,X^*,\Phi^\vee, \cdot^\vee)$, where $X^*=X^*(T)$ and $X_*=X_*(T)$ denote respectively the \emph{character} and \emph{cocharacter groups}, where $\Phi$ and $\Phi^\vee$ denote the \emph{root system} and \emph{coroot system}, and where $\cdot^\vee:\Phi\to\Phi^\vee:\alpha\mapsto\alpha^\vee$ is a bijection such that $\langle\alpha,\alpha^\vee\rangle=2$. For each $\alpha\in \Phi$, we denote by $u_\alpha:\G_a\simeq U_\alpha\subset G$ the corresponding \emph{root subgroup}.

With the choice of $T\subset B\subset G$, we have the subset $\Phi^+$ (resp.~$\Phi^-$) of positive (resp.~negative) roots in $\Phi$, i.e.~those $\alpha\in\Phi$ such that $U_\alpha\subset B$ (resp.~$U_\alpha\not\subset B$). Moreover, we have a subset $\Delta\subset\Phi^+$ of \emph{simple roots}. Every element in $\Phi^+=-\Phi^-$ can be expressed as a non-negative integer linear combination of elements in $\Delta$. Moreover, we have the \emph{root space decomposition} (also called the \emph{Cartan decomposition}) of the Lie algebra $\g$ of $G$ as $\g=\ft\oplus\bigoplus_{\alpha}\g_\alpha$, where $\ft=\mathrm{Lie}(T)$ is the Cartan sub-Lie algebra and $\g_\alpha$ is the (one-dimensional) $\alpha$-root space for each $\alpha\in \Phi$.

There exists a unique almost simple reductive group $\check{G}$  over $\C$ with root datum $(X^*,\Phi^\vee,X_*,\Phi, (\cdot^\vee)^{-1})$, up to isomorphism, called the \textit{dual group} of $G$. The corresponding maximal torus and Borel subgroup of $G$ are denoted by $\check{T}\subset \check{B}\subset \check{G}$. In particular, $X^*(\check {T})=X_*(T)$ and $X_*(\check{T})=X^*(T)$.

\paragraph{Acknowledgments.}
We are grateful to Claude Sabbah for answering our questions about various points in the theory of twistor modules.  We also thank Hao Sun, Daxin Xu, and Lingfei Yi for valuable discussions.

The first author was supported by the European Research Council (ERC) under the European Union’s Horizon 2020 research and innovation program (grant agreement no.\,101020009, project TameHodge). The second and third authors acknowledge support from the DFG-project SE 1114/6-1.

\section{Frenkel--Gross connections} \label{sec:FGconn}
In this section, we recall the definition of Frenkel--Gross connections in \cref{sec:FG-intro}, and discuss their basic properties such as functoriality in \cref{sec:functoriality} and their (exponentially) geometrical interpretation in \cref{sec:HNY-construction}. We cite Zhu's theorem (\cref{thm:zhu}) that expresses the Frenkel-Gross connection as a direct image of an exponentially twisted intersection complex. To finish, we will study pull-backs of Frenkel--Gross connections under cyclic covers in \cref{sec:homogeneity}. Using some specific homogeneity properties of the function by which the exponential twist is made, we give in \cref{prop:HNY-D-module} a variant of Zhu's theorem that will be crucial in the study of inverse Fourier transformations of Frenkel--Gross connections in \cref{sec:inverse-fourier}.

\subsection{Definition of Frenkel--Gross connections}\label{sec:FG-intro}
Let $\check{G}$ be an almost simple group over $\C$. As before, we fix a Borel subgroup $\check B\subset \check G$ and a maximal torus $\check{T}\subset \check{G}$. For each simple root $\alpha_i^\vee$ of $\check G$, we denote by $X_{-\alpha_i^\vee}$ a basis vector in the root subspace of $\check{\g}=\mathrm{Lie}(\check{G})$. Let $E=X_{\theta^\vee}$ be a basis vector in the root subspace of $\check{\g}$ corresponding to the maximal root $\theta^\vee$. Let $N=\sum_{\alpha^\vee\in \Delta^\vee} X_{-\alpha^\vee}$ be the corresponding principal nilpotent element. Then, the \emph{Frenkel--Gross connection}, denoted by $\nabla_{\check{G}}$, is a trivial $G$-bundle whose connection is defined by
\begin{equation}\label{eq:FG-conn}
    \nabla_{\check{G}}=\dif+\tfrac1t(N+Et)\, \dif t
\end{equation}
(Notice that the variable $t$ is called $z$ in \cite{FrenkelGross}, but we are following a different convention here to be consistent with the notation used in the theory of mixed twistor modules, especially with the convention we employ in \cref{sec:RMod}). Alternatively, we can view the Frenkel--Gross connection as a tensor functor
\begin{equation}\label{eq:functoriality}
    \nabla_{\check G}\colon \mathrm{Rep}(\check G)\to \mathrm{Conn}(\G_{m,t})
\end{equation}
from the category of finite-dimensional representations of $\check G$ to the category of flat connections on $\G_m$.

It was shown in \cite{FrenkelGross} that $\nabla_{\check G}$ has a regular singularity at $0$ with principal unipotent monodromy in the sense of \cite{Konstant}, with $N$ being the nilpotent part of the monodromy at $0$. Moreover, it has an irregular singularity at $\infty$ of slope $\frac{1}{h}$, where $h=h({\check G})$ is the Coxeter number of ${\check G}$. For the reason of slopes at $\infty$, we also work with
\begin{equation}\label{eq:FG-pull-back}
    \widetilde{\nabla}_{\check{G}}(V):=[h]^+\nabla_{\check G}(V)
\end{equation}
(where $[h]:\G_{m,t}\to\G_{m,t}$ denotes the $h$-th power map),
which has only $0$ or $1$ as the slopes at $\infty$, and can also be viewed as a tensor functor.

\begin{example}\label{exe:An}
\begin{itemize}
    \item  If $\check{G}=\SL_{n+1}$ and $V=\C^{n+1}$ the standard representation of $\SL_{n+1}$, then we recover the Bessel connection in \eqref{eq:Bessel}.
    \item  If $\check{G}=\rr{SO}_{2n+1}$ and $V=\bb{C}^{2n+1}$ the standard representation of $\rr{SO}_{2n+1}$, then $\nabla_{\rr{SO}_{2n+1}}(V)$ is isomorphic to the hypergeometric connection associated with the differential equation
    \[
        (t\partial_t)^{2n+1}-2t(t\partial_t+\tfrac12)=0,
    \]
    see  \cite[\S 6.3]{FrenkelGross}.
\end{itemize}
\end{example}

\subsection{Functoriality of Frenkel--Gross connections}\label{sec:functoriality}

The differential Galois groups of $\nabla_{\check{G}}$ are calculated in \cite[Cor.~9 \& 10]{FrenkelGross}. If $\check{G}$ is an almost simple group of a type listed on the left-hand side of the table below, then the differential Galois group $G_{\mathrm{gal}}$ is an almost simple group of the type listed on the right-hand side of the same row.

\begin{equation}\label{eq:monodromy-groups}
    \begin{tabular}{@{}l |cccccccc}
                        $\check{G}$ & $G_{\rr{gal}}$ \\
            \hline
            $A_{2n}$&$A_{2n}$\\
            $A_{2n-1}$, $C_n$& $C_n$\\
            $B_n$, $D_{n+1}(n\geq 4)$&$B_n$\\
            $E_7$ & $E_7$\\
            $E_8$ & $E_8$\\
            $E_6$, $F_4$ & $F_4$\\
            $B_3$, $D_4$, $G_2$ & $G_2$ \\
                \end{tabular}
\end{equation}

Let $\check{G}'\subset \check{G}$ be two almost simple groups that appear in the same line in the left column of \eqref{eq:monodromy-groups}. Then $\nabla_{\check{G}'}$ and $\nabla_{\check{G}}$ have the same differential Galois groups. So we can ask if $\nabla_{\check{G}'}$ is induced by $\nabla_{\check{G}}$. By \cite[Thm. 4.3.3 \& Prop. 5.13]{XZ22},  we have the following result:
\begin{theorem}
    For a choice of $N'$ and $E'$ in $\check{\mathfrak{g}}'$, there exist $N$ and $E$ in $\check{\mathfrak{g}}$ which coincide with $N'$ and $E'$ respectively under the inclusion $\check{\mathfrak{g}}'\subset \check{\mathfrak{g}}$. Moreover, $\nabla_{\check{G}}$ is the push-out of $\nabla_{\check{G}'}$, i.e.~the diagram
\begin{equation}\label{eq:pushout}
    \begin{tikzcd}
        \rr{Rep}(\check{G}) \ar[rr, "\text{restriction}"] \ar[rd,"\nabla_{\check{G}}"]& &\rr{Rep}(\check{G}') \ar[ld,"\nabla_{\check{G}'}"']\\
        & \mathrm{Conn}(\bb{G}_{m,t}) &
    \end{tikzcd}
\end{equation}
is commutative.
\end{theorem}
\begin{example}\label{exe:functoriality}
    \begin{enumerate}
        \item When $\check{G}'$ and $\check{G}$ are of the same type, they have the same Lie algebra. We have the `trivial functoriality' in this case. In other words, we have $\check{\mathfrak{g}}=\check{\mathfrak{g}}'$ and take $N=N', E=E'$. Then for each representation $V$ of $\check{G}$, we have $\nabla_{\check{G}}(V)=\nabla_{\check{G}'}(V).$
        \item When $\check{G}'=\rr{SO}_{2n+1}$ and $\check{G}=\rr{SO}_{2n+2}$, since the standard representation of $\mathrm{SO}_{2n+2}$ restricted as a representation of $\mathrm{SO}_{2n+1}$ splits into a direct sum of the standard representation of $\mathrm{SO}_{2n+1}$ and a trivial representation, we have
        \begin{equation}\label{eq:So2n+1-so2n+2}
            \nabla_{\mathrm{SO}_{2n+2}}(\dC^{2n+2})=\nabla_{\mathrm{SO}_{2n+1}}(\dC^{2n+1})\oplus \cO_{\bb{G}_{m}}.
        \end{equation}
        \item When $\check{G}'=G_2$ and $\check{G}=\rr{SO}_7$, the standard representation $V=\bb{C}^7$ of $\rr{SO}_7$ remains irreducible when viewed as a representation of $G_2$. In this case, we have
        \begin{equation}
            \nabla_{G_2}(V)=\nabla_{\rr{SO}_{7}}(V),
        \end{equation}
        which is the connection associated with the hypergeometric equation in \cref{exe:An}, see also \cite[\S6.4]{FrenkelGross} for the direct calculation in the case of $G_2$.
        \item When $\check{G}'=F_4$ and $\check{G}=E_6$, we have a decomposition $\mathfrak{e}_6=\mathfrak{f}_4\oplus \tilde{V}$, where $\tilde{V}$ is the $26$-dimensional         representation of $F_4$. Note that $\dim(E_6)=78$ and $\dim(F_4)=52$. Hence, the (nontrivial) $27$-dimensional minuscule $E_6$-representation $V$ restricts to a nontrivial representation of $F_4$, which therefore has to be $\bb{C}\oplus \tilde{V}$ (\textit{cf.} \cite[Lem. 14.4]{adams1998}). By the functoriality \eqref{eq:pushout}, we have
        \begin{equation}\label{eq:decomposition-f4-e6}
             \nabla_{E_6}(V)=\nabla_{F_4}(\widetilde V)\oplus \cO_{\bb{G}_{m}}
        \end{equation}
        for compatible choices of $N$ and $E$ on both sides.
    \end{enumerate}
\end{example}

\subsection{Heinloth--Ngo--Yun's Kloosterman connections}\label{sec:HNY-construction}
\label{sec:HNYconstruction}

Heinloth--Ng\^o--Yun constructed ($\ell$-adic) Kloosterman sheaves for reductive groups in \cite{HeinlothNgoYun}. Applying their construction to $\mathscr{D}$-modules, we also get some $\check{G}$-connections on $\bb{G}_m$.

For each almost simple group $G$, we fix a maximal torus and a Borel subgroup $T\subset B\subset G$. Recall that the \textit{loop group} $LG$ (resp. the \textit{positive loop group} $L^+G$) is the fppf sheaf on the category of $k$-algebras, defined by
    \[
    R\mapsto G(R((t)) \quad (\text{resp. } R\mapsto G(R[[t]])).
    \]
The \textit{affine Grassmannian} $\rr{Gr}=\rr{Gr}_G$ is the fppf-quotient $LG/L^+G$ \cite[Prop.~1.3.6]{Zhu17Grassmannian}, which is represented by an ind-scheme.
The group $L^+G$ has a natural action on $\rr{Gr}$, and the orbits are indexed by dominant cocharacters $\lambda\in \Phi^{\vee+}$. This orbit and its closure are denoted by $\rr{Gr}_{\lambda}$ and $\rr{Gr}_{\leq \lambda}$ respectively \cite[Prop.~2.1.5]{Zhu17Grassmannian}.
$\rr{Gr}_{\leq \lambda}$ is a finite-dimensional algebraic variety,
with $\dim(\rr{Gr}_{\leq \lambda})=\langle2 \rho,\lambda\rangle$, and
$\rr{Gr}_{\lambda}$ is an open dense subvariety of it.
We denote by $\rr{IC}_{\lambda}$ the intersection complex on $\rr{Gr}_{\leq \lambda}$.
For each irreducible representation $V$ of the dual group $\check{G}$ of highest weight $\lambda$, we also write $\rr{IC}_V=\rr{IC}_{\lambda}$.

The Beilinson--Drinfeld Grassmanian,
$\pi\colon \rr{GR}\to \bb{G}_{m}$, over $\bb{G}_{m}$ is a global version of affine Grassmannians, see for example \cite[Def. 3.11]{Zhu17Grassmannian}. We denote by $\bb{G}_{m}^{\rr{rot}}$ the torus $\bb{G}_{m}$ that acts on $\bb{G}_{m}$ by multiplication. The specified action of $\bb{G}_{m}^{\rr{rot}}$ on $\bb{G}_m$  induces an action of $\bb{G}_m^{\rr{rot}}$ on $\rr{GR}$ such that $\pi$ is $\bb{G}_{m}^{\rr{rot}}$-equivariant. Using this action, we can trivialize the fibration $\pi$ as
\[
     \rr{GR}\simeq \rr{Gr}\times \bb{G}_m \to \bb{G}_{m},
\]
where $\rr{Gr}$ is the affine Grassmannian for $G$, see \cite[(5.11)]{HeinlothNgoYun}.

Let $\rr{Bun}_{G(1,2)}$ be the moduli stack of $G$-bundles on $\bb{G}_m$ with some level structures at $0$ and $\infty$ determined by level groups $I(1)$ and $I(2)$, as defined in \cite[\S1.2]{HeinlothNgoYun}. Then there is a projection $\mathrm{pr}\colon \rr{GR}\to \rr{Bun}_{G(1,2)}$, and we have the following diagram in \cite[\S5.2]{HeinlothNgoYun} or \cite[(2.2)]{Yun15} as follows:
\begin{equation}\label{eq:HNY-simplified}
    \begin{tikzcd}
     &&\rr{GR}^\circ \ar[r,hook]\ar[lld,bend right=15,"F"']\ar[ld,"f"]\arrow[rrd, bend left=55,"\pi^\circ"] &\rr{GR} \ar[rd,"\pi"'] \ar[ld,"\rr{pr}"]  & \\
       \bb{A}^1 &\bb{G}_a^{r+1}\ar[r,hook,"\phi"']\ar[l,"\rr{sum}"] &\rr{Bun}_{G(1,2)}&&\bb{G}_{m}
\end{tikzcd}
\end{equation}
where $\phi\colon \bb{G}_a^{r+1}\simeq T\times I(1)/I(2)\hookrightarrow \rr{Bun}_{G(1,2)}$ is the inclusion of the big open cell \cite[Cor.\,1.3\,(4)]{HeinlothNgoYun}; the open sub-ind-scheme $\rr{GR}^\circ\subset \rr{GR}$ is the inverse image of the big cell.

Similar to the Beilinson--Drinfeld Grassmanian, we have $\rr{GR}^\circ\simeq \rr{Gr}^\circ\times \bb{G}_m$. Let $\pi^\circ\colon \rr{GR}^\circ \to \bb{G}_m$ and $\pi_{\rr{Gr}}\colon \rr{GR}\to \rr{Gr}$ be the corresponding projections, and we write $\rr{IC}_{V,\rr{GR}^\circ}:=\pi_{\rr{Gr}}^+\rr{IC}_V|_{\rr{Gr}^\circ}$.
Then the two isomorphic complexes of $\mathscr{D}$-modules
    \begin{equation}\label{eq:!vs*}
    \pi^\circ_\dagger(\cc{E}^F\otimes \rr{IC}_{V,\rr{GR}^\circ})\simeq \pi^\circ_+(\cc{E}^F\otimes \rr{IC}_{V,\rr{GR}^\circ})
    \end{equation}
are connections over $\bb{G}_{m}$, denoted by $\Kl_{\check{G}}(V)$, where $\mathcal{E}^F=(\mathcal{O}_{\rr{GR}^\circ},\rr{d}+\rr{d}F)$ is the exponential $\scr{D}$-module. Furthermore, they can be upgraded to a tensor functor
\[
    \Kl_{\check{G}}\colon \rr{Rep}(\check{G})\to \rr{Conn}(\bb{G}_{m}),
\]
see \cite[Thm. 1(1)]{HeinlothNgoYun} or \cite[Thm. 2.2.1]{Yun15}.

Heinloth--Ngô--Yun conjectured that these $\scr{D}$-modules should coincide with Frenkel--Gross's connections \cite[Conj. 2.16]{HeinlothNgoYun}, and this conjecture has since been proven:
\begin{theorem}[{\cite[\S6]{Zhu17}}]\label{thm:zhu}
    The tensor functor $\Kl_{\check{G}}$ and the Frenkel--Gross connection $\nabla_{\check{G}}$ are isomorphic.
\end{theorem}

\subsection{Homogeneity properties} \label{sec:homogeneity}
Let $G$ be an almost simple group of adjoint type, so the half sum of positive coroots $\rho^\vee$ is a cocharacter, and recall that we have fixed a maximal torus $T\subset G$. In \cite[\S 2.6.1]{Yun15}, Yun constructed an ind-scheme $\mathfrak{S}=\mathfrak{S}_1\times \bb{G}_{m}$ isomorphic to $\rr{GR}^\circ=\rr{Gr}^\circ\times \bb{G}_{m}$, compatible with the $\bb{G}_{m}^{\rr{rot}}$-action. In \cite[\S2.6.4]{Yun15}, he also defined an embedding of a one-dimensional subtorus $\bb{G}_m$ of $T\times \bb{G}_{m}^{\rr{rot}}$, labeled $\bb{G}_{m}^{(\rho^\vee,h)}$, by
    $$\bb{G}_{m}^{(\rho^\vee,h)}\ni t \mapsto (\rho^\vee(t),t^h)\in T\times \bb{G}_{m}^{\rr{rot}}.$$
Hence, the torus $\bb{G}_{m}^{(\rho^\vee,h)}$ inherits an action on $\rr{GR}^\circ=\rr{Gr}^\circ\times \bb{G}_{m}$. Moreover, $F$ is $\bb{G}_{m}^{(\rho^\vee,h)}$-equivariant and $\bb{G}_{m}^{(\rho^\vee,h)}$ acts on both $\bb{G}_a^{r+1}$ and $\bb{A}^1$ by the multiplication on each factor.

Recall that $[h]\colon \bb{G}_{m,t}\to \bb{G}_{m,t}$ denotes the $h$th power morphism. Let $\rr{act}\colon  \rr{Gr}^\circ\times \bb{G}_{m}\to \rr{Gr}^\circ$ be the action of $\bb{G}_{m}^{(\rho^\vee,h)}$ on $\rr{Gr}^\circ$ and $\tilde{F}=F\circ (\rr{act}\times [h])$. These morphisms are illustrated in the following diagram:
\begin{equation}\label{eq:HNY-diagram}
\begin{tikzcd}
     &\rr{Gr}^\circ \times \bb{G}_{m,t} \ar[ldd,bend right,"\tilde{F}"'] \ar[d,"\rr{act}\times {[h]}"] \ar[rd,"\rr{act}\times \rr{id}",bend left] &&\\
     &\rr{Gr}^\circ \times \bb{G}_{m} \ar[rd,"\pi^\circ"] \ar[ld,"F"'] & \rr{Gr}^\circ \times \bb{G}_{m,t} \ar[rd,"\pi^\circ"]\ar[l,"\rr{id}\times {[h]}"'] & \\
        \bb{A}^1 &&\bb{G}_{m} & \bb{G}_{m,t} \ar[l,"{[h]}"'].
\end{tikzcd}.
\end{equation}

Then we deduce that
    \begin{align*}
                        \tilde{F}(x,t)&            =F(\rr{act}(x,t),t^h) &&\text{definition of $\tilde{F}$}\\
            &=F(t\cdot (x,1)) &&\text{definition of the action of $\bb{G}_{m}^{(\rho^\vee,h)}$}\\
            &=tF(x,1) &&\text{$F$ is $\bb{G}_{m}^{(\rho^\vee,h)}$-equivariant.}
    \end{align*}
Notice that $\pi^\circ\circ (\rr{act}\times \rr{id})=\pi^\circ$, we have
\begin{equation}\label{eq:HNY-D-module}
    \begin{split}
        [h]^+\rr{Kl}_{\check{G}}(V)
    &= [h]^+\pi^\circ_+(\cc{E}^F\otimes \rr{IC}_{V,\rr{GR}^\circ})
        \cong \pi^\circ_+(\cc{E}^{F\circ(\rr{id}\times [h])}\otimes \rr{IC}_{V,\rr{GR}^\circ})\\
    &\cong \pi^\circ_+(\rr{act}\times \rr{id})_+(\rr{act}\times \rr{id})^+(\cc{E}^{F\circ(\rr{id}\times [h])}\otimes \rr{IC}_{V,\rr{GR}^\circ})\\
    &\cong \pi^\circ_+(\cc{E}^{\tilde{F}}\otimes \rr{IC}_{V,\rr{GR}^\circ})
        \cong {\pi}^\circ_+(\cc{E}^{t\cdot F(x,1)}\otimes \rr{IC}_{V,\rr{GR}^\circ}),
    \end{split}
\end{equation}
where in the third isomorphism, we used the fact that
$\rr{act}\times \rr{id}$ is an isomorphism.

Combined with \cref{thm:zhu}, we have the following proposition.
\begin{proposition}\label{prop:HNY-D-module}
    The connection $\widetilde\nabla_{\check{G}}(V)=[h]^+\nabla_{\check{G}}(V)$ is of the form
    \[
        {\pi}^\circ_+(\cc{E}^{t\cdot g}\otimes \rr{IC}_{V,\rr{GR}^\circ}),
    \]
    where $g=F|_{\rr{Gr}^{\circ}\times\{1\}}\colon \mathrm{Gr}^\circ \to \bb{A}^1:x\mapsto F(x,1)$, $t\cdot g$ is viewed as a regular function on $\mathrm{Gr}^\circ\times \bb{G}_{m,t}$, and $\pi^\circ$ is the projection from  $\mathrm{Gr}^\circ\times \bb{G}_{m,t}$ to $\bb{G}_{m,t}$.
\end{proposition}

\section{The irregular Hodge numbers of Frenkel--Gross connections} \label{sec:IrregHodge}
In this section, the aim is to prove \cref{thm:main}. We first recall necessary properties of exponential mixed Hodge modules (or rescalable integrable mixed twistor $\scr{D}$-modules) in \cref{sec:MTM}. Then we will equip Frenkel--Gross connections with structures of exponential mixed Hodge modules in \cref{sec:EMHM}. In \cref{sec:inverse-fourier}, we will calculate the Hodge numbers of inverse Fourier transformations of pull-backs of Frenkel--Gross connections, which can be calculated from the unipotent monodromy. In \cref{sec:proof}, we use the functorial properties of Frenkel--Gross connections and a filtered stationary phase formula due to Sabbah-Yu to conclude the proof of \cref{thm:main}.

\subsection{Recollections on rescalable integrable mixed twistor \texorpdfstring{$\scr{D}$}{D}-modules}
\label{sec:MTM}

Recall from the introduction that we have for any complex algebraic variety $M$ the abelian category $\mathrm{MHM}(M,\dR)$ of \textit{mixed Hodge modules} on $M$ with coefficients in $\dR$, which is the category of $\dR$-polarized mixed Hodge structures if $M$ is a point. Recall also that we have a forgetful functor $\mathrm{MHM}(M,\dR) \to \scr{D}_M$ (and an induced functor between the respective derived categories) compatible with the six functor formalism on both sides.

There is a related category $\mathrm{MTM}^{\rr{int}}(M,\dR)$ of \textit{integrable mixed twistor $\scr{D}_M$-modules} when $M$ is smooth, which is built in a parallel way (with major new challenges though) to the theory of mixed Hodge modules but uses $\scr{R}_M^\rr{int}=\scr{D}_M[z]\langle z^2\partial_z \rangle$-modules instead of $\scr{D}_M$-modules. The underlying $\scr{R}$-modules $\scr{M}$ of objects in $\mathrm{MTM}^{\rr{int}}(M)$  give rise to (not necessarily regular) holonomic $\scr{D}_M$-modules through the functor $\Xi_{\rr{dR}}(\scr{M})=\scr{M}/(z-1)\scr{M}$. The induced functor $\Xi_{\rr{dR}}\colon D^b\rr{MTM}^{\rr{int}}(M)\to D^b_h(\scr{D}_M)$ is also compatible with the six functor formalism and nearby cycle functors \cite[Paragraph after Prop. 14.1.24]{Mochizuki-MTM}. Furthermore, there is a fully faithful exact functor
\begin{equation}\label{eq:Rees-construction}
    v\colon D^b(\mathrm{MHM}(M)) \to D^b\mathrm{MTM}^{\rr{int}}(M)
\end{equation}
that is built via the Rees construction and which is compatible with the six-functor formalism, as detailed in \cite[Prop.\,13.5.5 \& Prop.\,14.3.29]{Mochizuki-MTM}.

We are interested in a category of ``irregular Hodge modules'' that lies in between $\rr{MHM}(M,\dR)$ and $\rr{MTM}^{\rr{int}}(M)$. It was first defined and studied by Sabbah \cite{Sabbah-Irreg-book}. We will rather work with two related categories (studied by Mochizuki in \cite{Mochizuki}): namely, the category of rescalable integrable twistor $\scr{D}$-modules and the category of exponential mixed Hodge modules.

\begin{definition}\label{def:ExpHM}
    The category $\rr{EMHM}(M)$ of \textit{exponential mixed Hodge modules} on $M$ is defined as the full subcategory of $\mathrm{MHM}(\bb{A}^1_\theta\times M,\dR)$ whose objects $N^{\mathrm{H}}$ satisfy $\pi_{M*}N^{\mathrm{H}}=0$ for the projection $\pi_M\colon \bb{A}^1_\theta\times M\to M$. When $M$ is a point, we recover the category of exponential mixed Hodge structures, as defined in \cite{Kontsevich-Soibelman}.
\end{definition}

The category $\EMHM(M)$ can also be defined as the essential image of a projector $P_*\colon \mathrm{MHM}(\bb{A}^1_\theta\times M,\dR)\to \mathrm{MHM}(\bb{A}^1_\theta\times M,\dR)$, which is defined in \cite[(155)]{Mochizuki},
see also the discussion right after it.

Now we explain how to embed $\rr{EMHM}(M)$ as a subcategory of $\rr{MTM}^{\rr{int}}(M)$ and to attach an irregular Hodge filtration to its objects.  Mochizuki constructed a functor
\begin{equation}\label{eq:rescalable-functor}
    u\colon \rr{MTM} ^{\rr{int}}(\bb{A}^1_\tau \times M) \to \rr{MTM} ^{\rr{int}}(M)
\end{equation}
in \cite[\S11.4.4]{Mochizuki} by sending $\scr{M}$ to an object $u(\scr{M})$ such that $\iota_!u(\scr{M})=\psi^{(1)}_{\tau-1}\scr{M}$.
\begin{definition}[{\cite[Def. 11.38]{Mochizuki}}]
    The essential image of $u$ is a full subcategory of $\rr{MTM} ^{\rr{int}}(M)$, called \textit{rescalable mixed twistor $\scr{D}_Y$-modules} and denoted by $\rr{MTM}_{\rr{resc}}^{\rr{int}}(M)$.
\end{definition}

Let $\mathbf{FT}_M \colon \rr{MTM} ^{\rr{int}}(\bb{A}^1_\theta \times M) \to \rr{MTM} ^{\rr{int}}(\bb{A}^1_\tau \times M) $ be the partial Fourier transform \cite[\S10.2]{Mochizuki} of integrable twistor modules relative to $M$ where $\tau$ denotes the coordinate of the affine line dual to $\bb{A}^1_\theta$. Mochizuki proved the following theorem:

\begin{theorem}[{\cite[Thm.\,11.45]{Mochizuki}}]\label{thm:Mochizuki}
    The functor
    \begin{equation}
        B\colon \rr{MHM}(\bb{A}^1_{\theta}\times M,\dR) \xrightarrow{u\circ \mathbf{FT}_M\circ v} \rr{MTM}^{\rr{int}}(M),
    \end{equation}
    induces an equivalence of categories of $B_*\colon \rr{EMHM}(M)\xrightarrow{\sim}\rr{MTM}_{\rr{resc}}^{\rr{int}}(M)$.
\end{theorem}

Moreover, Mochizuki showed that the functor $B$ is compatible with the six-functor formalism, nearby and vanishing cycle functors, and duality in \cite[Props.~11.46-11.48]{Mochizuki}. In particular, these operators preserve rescalable integrable mixed twistor $\scr{D}$-modules. Let us also notice that the category $\MHM(\bb{A}^1_{\theta}\times M,\dR)$ exists even if $M$ is a possibly singular algebraic variety. Furthermore, the push-forward and the pull-back for maps between possibly singular varieties preserve the subcategory of exponential mixed Hodge modules. Therefore, we have the six-functor formalism in the category of exponential mixed Hodge modules (and consequently in the category of rescalable integrable mixed twistor $\cD$-modules) including push-forward and pull-back for maps between singular spaces.

As explained above, the \textit{de Rham} functor $\Xi_{\rr{dR}}:
\MTM^{\rr{int}}_{\rr{resc}}(M)\subset \MTM(M) \rightarrow \textup{Mod}_h(\cD_Y)$ sends a rescalable integrable twistor $\scr{D}$-module $\cM$ to the holonomic $\scr{D}$-module $\cM/(z-1)\cM$. For any $\cT\in\MTM^{\rr{int}}_{\rr{resc}}(M)$ one can then equip $\Xi_{\rr{dR}}(\cT)$ with an irregular Hodge filtration $F_{\rr{irr}}^\bullet$, see \cite[Thm.~0.3]{Sabbah-Irreg-book} and \cite[Cor.~1.6]{Mochizuki}. In particular, we denote by $h^p_{\rr{irr}}$ the irregular Hodge numbers, i.e. the (generic) ranks of the graded pieces of $F_{\rr{irr}}^\bullet \Xi_{\rr{dR}}(\cT)$.

To end this section, we collect a unicity result, which will be needed later.

\begin{lemma}\label{lem:irreducible}
Let $M$ be a smooth quasi-projective variety, and let $N$ be a semi-simple holonomic $\cD_Y$-module. Then there is up to isomorphism at most one object $\cT\in \rr{MTM}^{\rr{int}}_{\rr{resc}}(M)$ of weight $0$ such that $\Xi_{\rr{dR}}(\cT)\cong N$.
\end{lemma}

\begin{proof}
We choose a compactification $j_M:M\hookrightarrow \overline{M}$, where $\overline{M}$ is smooth and projective. Then the intermediate extension $j_{M,\dag+}N$ is a semi-simple holonomic $\cD_{\overline{M}}$-module. Now, suppose that we have $\cT_1,\cT_2\in\rr{MTM}^{\rr{int}}_{\rr{resc}}(M)$ such that $\Xi_{\rr{dR}} (\cT_1) \cong \Xi_{\rr{dR}} (\cT_2)$.  According to \cite[Thm.~14.3.16 \& Prop.~14.3.17]{Mochizuki-MTM}, for $i\in\{1,2\}$, we have the objects $j_{M,*} \cT_i$ and $j_{M,!} \cT_i$ in $\MTM(\overline{M})$ and by loc.cit.~Proposition 14.3.18 there is a morphism $j_{M,!}\cT_i\rightarrow j_{M,*}\cT_i$. Denote by $j_{M,!*}\cT_i$ the respective images in $\MTM(\overline{M})$. Since taking the de Rham functor $\Xi_{\textup{dR}}$ commutes with taking direct images (with and without proper support) of integrable mixed twistor modules resp.~of holonomic $\cD$-modules, we see that $\Xi_{\rr{dR}}(j_{M,!*}\cT_i)=j_{M,\dag+}N$ for both $i=1,2$. Now by \cite[Thm. 1.4.4]{Mochizuki-Wild}, there is a unique (wild and pure) analytic twistor module $\overline{\cT}$ on $\overline{M}$ such that $\Xi_{\rr{dR}}(\overline{\cT})=(j_{M,\dag+}N)^{an}$. Hence, we obtain an isomorphism
$$
\phi:(j_{M,!*}\cT_1)^{an} \longrightarrow (j_{M,!*}\cT_2)^{an},
$$
but since $\overline{M}$ is projective, the categories $\MTM(M)^{an}$ and $\MTM(\overline{M})$ are equivalent (this follows, e.g., from \cite[Lem. 14.1.2]{Mochizuki-MTM}). Therefore we obtain an isomorphism
$$
\phi^{\rr{alg}}:j_{M,!*}\cT_1 \longrightarrow j_{M,!*}\cT_2,
$$
which restricts to an isomorphism
$$
\phi^{\rr{alg}}_{|M}:\cT_1 \longrightarrow \cT_2,
$$
of (algebraic) mixed twistor modules on $M$. If we now consider $\cT_i$
as integrable twistor modules on $M$, then by \cite[Rem. 1.39]{Sabbah-Irreg-book}, $\phi^{\rr{alg}}_{|M}$ is also an isomorphism in $\MTM^{\rr{int}}(M)$.
Since $\MTM_{\rr{resc}}^{\rr{int}}(M)$ is a a full subcategory of $\MTM^{\rr{int}}(M)$, $\phi_{|M}$ is an isomorphism of (algebraic) rescalable integrable mixed twistor $\scr{D}$-modules, as required.
\end{proof}

\subsection{Rescalable integrable mixed twistor \texorpdfstring{$\scr{D}$}{D}-modules attached to \texorpdfstring{$\nabla_{\check{G}}(V)$}{}.}\label{sec:EMHM}

In this section, we construct objects in $\MTM^{\rr{int}}_{\rr{resc}}(\bb{G}_{m,t})$, which will be sent to Frenkel--Gross connections under the de Rham functors.

For a regular function $f\colon U\to \bb{A}^1$, there exists an integrable twistor $\scr{D}$-modules $\cc{T}^{f/z}\in \rr{MTM}^{\rr{int}}(U)$, which is sent to the exponential $\scr{D}_U$-module $\cc{E}^f=(\cc{O}_U,\rr{d}+\rr{d}f)$ under the de Rham functor. See \cite[Discussion before Thm.~0.2]{Sabbah-Irreg-book} for details. One can verify that it is rescalable, as it is the image of $\cc{T}^{\tau \cdot f/z}$ under the functor $u$ mentioned above in \eqref{eq:rescalable-functor}. (Note that $\cT^{\tau\cdot f/z}$ is a smooth $\cR_{\dA^1_\tau\times U}$-module, so that $\psi^{(1)}_{\tau-1} \cT^{\tau\cdot f/z}=\cT^{f/z}$.)

\paragraph{The case of \texorpdfstring{${\nabla}_{\check{G}}(V)$}{}}

Let $V$ be an irreducible $\check{G}$-representation of the highest weight $\lambda$, then $\rr{Gr}_{\leq \lambda}$ has dimension $d_V:=\langle 2\rho,\lambda\rangle$. There exists a mixed Hodge module with underlying $\scr{D}$-module $\rr{IC}_{V,\rr{GR}^\circ}$, denoted by $\rr{IC}^{\rr{H}}_{V,\rr{GR}^\circ}$. In particular, after tensoring with a suitable Tate's shift, the mixed Hodge module $\rr{IC}^{\rr{H}}_{V,\rr{GR}^\circ}(\frac{d_V+1}{2})$ has weight $0$. Taking notation from \eqref{eq:HNY-simplified}, we define
    \begin{equation}\label{eq:EMHM-nabla}
        \nabla_{\check{G}}^{\rr{H}}(V):= \pi^\circ_*\Bigl(\cc{T}^{F/z}\otimes v(\rr{IC}_{V,\rr{GR}^\circ})\left(\tfrac{d_v+1}{2}\right)\Bigr)
    \end{equation}
where $v$ is the Rees functor \eqref{eq:Rees-construction}.

\begin{remark}
    Although $\rr{GR}^\circ$ is infinite dimensional, the support of $\rr{IC}_{V,\rr{GR}^\circ}^{\rr{H}}$ is the finite dimensional variety $\rr{GR}^\circ \cap \rr{GR}_{\leq \lambda}$. Therefore, in \eqref{eq:EMHM-nabla} and similarly in the sequel of this paper we understand the morphisms $F$ and $\pi^\circ$ to be restricted to the support of $\rr{IC}_{V,\rr{GR}^\circ}^{\rr{H}}$. See also the discussion after \cref{thm:Mochizuki} about the six-functor formalism for rescalable integrable mixed twistor $\cD$-modules.
\end{remark}

By the compatibility of the six-functor formalism with the de Rham functor $\Xi_{\rr{dR}}$ we conclude from \eqref{eq:!vs*}
\[
    \Xi_{\rr{dR}}\Bigl( {\nabla}_{\check{G}}^{\rr{H}}(V)\Bigr)= \pi^\circ_+(\cc{E}^{F}\otimes  \rr{IC}_{V,\rr{GR}^\circ})|_{\bb{G}_{m,t}}= \nabla_{\check{G}}(V).
\]
Hence, we call ${\nabla}_{\check{G}}^{\rr{H}}(V)$ the  rescalable integrable twistor $\scr{D}$-module attached to $\nabla_{\check{G}}(V)$.
This also explains the choice of notation: we follow the
convention from \cite{Saito1, SaitoMHM} to denote by $\nabla_{\check{G}}^{\rr{H}}(V)$ the Hodge/twistor object, which has $\nabla_{\check{G}}(V)$ as underlying $\cD$-module (and a similar convention applies to the Hodge module $\rr{IC}^{\rr{H}}_{V,\rr{GR}^\circ}$ from above).

Moreover, in \eqref{eq:!vs*}, the forget support morphism $\pi^\circ_\dagger(\cc{E}^{F}\otimes  \rr{IC}_{V,\rr{GR}^\circ})|_{\bb{G}_{m,t}}\to \pi^\circ_+(\cc{E}^{F}\otimes  \rr{IC}_{V,\rr{GR}^\circ})|_{\bb{G}_{m,t}}$ is an isomorphism. As $\Xi_{\rr{dR}}$ is an exact functor, the forget support morphism
\[
\pi^\circ_!\Bigl(
     \cT^{F/z}\otimes  v(\rr{IC}_{V,\rr{GR}^\circ}) \Bigr)|_{\bb{G}_{m,t}}\to \pi^\circ_*\Bigl(
     \cT^{F/z}\otimes  v(\rr{IC}_{V,\rr{GR}^\circ}) \Bigr)|_{\bb{G}_{m,t}}
\]
is also an isomorphism. Hence, by the formalism of weights,  $ {\nabla}_{\check{G}}^{\rr{H}}(V)$ is a pure twistor module of weight $0$.

\paragraph{The case of \texorpdfstring{$\widetilde{\nabla}_{\check{G}}(V)$}{}.}

Let $[h]\colon \bb{G}_{m,t}\to \bb{G}_{m,t}$ be the $h$-th power map. We define
\begin{equation}\label{eq:tilde-nabla-V}
    \widetilde{\nabla}_{\check{G}}^{\rr{H}}(V):=[h]^*\nabla^{\rr{H}}_{\check{G}}(V) \in \rr{MTM}^{\rr{int}}_{\rr{resc}}(\bb{G}_{m,t}).
\end{equation}
By the compatibility of the six-functor formalism with the de Rham functor $\Xi_{\rr{dR}}$ we conclude that
\[
    \Xi_{\rr{dR}}\Bigl( \widetilde{\nabla}_{\check{G}}^{\rr{H}}(V)\Bigr)= [h]^+\nabla_{\check{G}}(V)=\widetilde\nabla_{\check{G}}(V).
\]
Hence, we call $ \widetilde{\nabla}_{\check{G}}^{\rr{H}}(V)$ the  rescalable integrable twistor $\scr{D}$-module attached to $\widetilde\nabla_{\check{G}}(V).$

When $\check{G}$ is simply-connected (i.e.~$G$ is adjoint), we have by homogeneity (see \cref{prop:HNY-D-module}) that
\begin{equation}
    \widetilde{\nabla}_{\check{G}}^{\rr{H}}(V)= \pi^\circ_*\left(
     \cT^{(t\cdot g)/z}\otimes  v(\rr{IC}_{V,\rr{GR}^\circ})\left(\tfrac{d_V+1}{2}\right)\right)|_{\bb{G}_{m,t}}.
\end{equation}

\paragraph{The case of \texorpdfstring{$\widetilde{\nabla}_{\check{G}}(V)\otimes \bb{L}_\chi$}{}.}

For a rational number $a\in \dQ\backslash \dZ$, let $\chi=\exp(2\pi i a)\neq 1$ and let $\bb{L}_\chi$ be the connection defined by $\bb{L}_{\chi}:=(\mathcal{O}_{\bb{G}_{m,t}},\mathrm{d}+a\frac{\mathrm{d}t}{t})$ (sometimes called Kummer module). We then have the $\cD_{\dA^1_t}$-module $\widetilde{\nabla}_{\check{G}}(V)_\chi:=j_+( \widetilde{\nabla}_{\check{G}}(V)\otimes \bb{L}_\chi)$.
On the other hand, $\bb{L}_\chi$ underlies a pure complex Hodge module of rank $1$ that we denote by $\bb{L}_\chi^{\rr{H}}$.

We then consider the rescalable integrable mixed twistor $\scr{D}$-module
    \begin{equation}\label{eq:EMHM-chi}
        \widetilde{\nabla}_{\check{G}}^{\rr{H}}(V)_{\chi}:=j_{*}\widetilde{\nabla}_{\check{G}}^{\rr{H}}(V)\otimes j_*v(\bb{L}_\chi^{\rr{H}})\in \rr{MTM}_{\rr{resc}}^{\rr{int}}(\bb{A}^1_t)
    \end{equation}
that is attached to $\widetilde{\nabla}_{\check{G}}(V)_\chi$ in the sense that $\Xi_\rr{dR}(\widetilde{\nabla}_{\check{G}}^{\rr{H}}(V)_{\chi})=\widetilde{\nabla}_{\check{G}}(V)_\chi$.

\paragraph{Tensor functor.}

Similar to the construction of Kloosterman $\scr{D}$-modules in \cite[\S 2.6]{HeinlothNgoYun}, the rescalable integrable mixed twistor $\scr{D}$-module $\nabla_{\check{G}}^{\rr{H}}(V)=\pi^\circ_*\bigl( \mathcal{T}^{  F/z}\otimes \mathrm{IC}_{V,\rr{GR^\circ}}^{\rr{H}}\bigr)$
can be obtained by imitating the construction of Hecke-eigensheaves of the Hecke operator $\rr{Hk}_V$ in \cite[\S2.3]{HeinlothNgoYun}.

More precisely, we can construct the Hecke operator for rescalable integrable mixed twistor $\scr{D}$-modules, by interpreting \cite[page 259, last displayed formula]{HeinlothNgoYun} within the category rescalable integrable mixed twistor $\scr{D}$-modules. For a representation $V$ of $\check{G}$, we denote the corresponding functor from $\rr{MTM}^{\rr{int}}_{\rr{resc}}(\mathrm{Bun}_{G(1,2)})$ to $\rr{MTM}^{\rr{int}}_{\rr{resc}}(\mathrm{Bun}_{G(1,2)}\times \bb{P}^1\backslash\{0,\infty\})$ by $\rr{Hk}^\rr{H}_V$.

Moreover, the rescalable integrable mixed twistor $\scr{D}$-module $A_{\phi,\chi}^{\mathrm{H}}$ can be constructed by replacing the Artin--Schreier sheaf $\scr{L}_{\psi(x)}$ with $\cc{T}^{x/z}$ in the definition of $A_{\phi,\chi}$ in \cite[\S2.2]{HeinlothNgoYun}.  By the compatibility of the six-functor formalism and the de Rham functor, we deduce for any representation $V$ of $\check{G}$ that
\[
    \Xi_{\rr{dR}}(\rr{Hk}^\rr{H}_V(A_{\phi,\chi}^{\mathrm{H}}))=\rr{Hk}_V(A_{\phi,\chi})=A_{\phi,\chi}\boxtimes \nabla_{\check{G}}(V)
\]
in $D^b(\scr{D}_{\mathrm{Bun}_{G(1,2)}\times \bb{P}^1\backslash\{0,\infty\}})$, where $G$ is the dual group of $\check{G}$, and $\mathrm{Bun}_{G(1,2)}$ is the moduli stack of $G$-bundles with some level structures.
By \cref{lem:irreducible}, semi-simple holonomic $\mathscr{D}$-modules can underlie at most one rescalable integrable mixed twistor $\scr{D}$-module of weight $0$. Hence, we deduce that
\[
    \rr{Hk}^\rr{H}_V(A_{\phi,\chi}^{\mathrm{H}})=A^{\mathrm{H}}_{\phi,\chi}\boxtimes \nabla^{\mathrm{H}}_{\check{G}}(V).
\]
As in \cite[\S2.3]{HeinlothNgoYun}, we can compose Hecke operators using the geometric Satake correspondence for perverse sheaves (resp.~regular holonomic $\cD$-modules). More precisely, for any two representations of $\check{G}$, the composition $\rr{Hk}_V\circ \rr{Hk}_W$ equals the functor $\rr{Hk}_{V\otimes W}$ because the convolution of $\mathrm{IC}_V$ with $\mathrm{IC}_W$ is isomorphic to $\mathrm{IC}_{V\otimes W}$ (in particular, the convolution of $\mathrm{IC}_V$ with $\mathrm{IC}_W$ is still semi-simple). Using again \cref{lem:irreducible} for these intersection complexes, we deduce that the convolution of $\mathrm{IC}^{\mathrm{H}}_V$ with $\mathrm{IC}^{\mathrm{H}}_W$ is isomorphic to $\mathrm{IC}^{\mathrm{H}}_{V\otimes W}$, which implies that $\rr{Hk}^\rr{H}_V\circ \rr{Hk}^\rr{H}_W=\rr{Hk}^\rr{H}_{V\otimes W}$. Hence, $\nabla_{\check{G}}^{\rr{H}}$ yields a tensor functor
\begin{equation}\label{eq:EMHM-tensor}
    \nabla_{\check{G}}^{\rr{H}}\colon \rr{Rep}_{\check{G}} \to \rr{MTM}^{\rr{int}}_{\rr{resc}}(\bb{G}_{m,t}).
\end{equation}

As explained in \cref{sec:EMHM}, each object in the category of rescalable mixed twistor modules is equipped with an irregular Hodge filtration in a functorial way. When restricted to smooth objects (that is, mixed twistor $\scr{D}$-modules whose image under the de Rham functor are smooth $\cD$-module, i.e. connections), this gives a tensor functor
$$
    \rr{MTM}^{\rr{int}}_{\rr{sm,resc}}(\bb{G}_{m,t}) \xrightarrow{\rr{IrrFil}}
\rr{Fil}_{\bb{G}_{m,t}}
$$
to the category of vector bundles on $\G_{m,t}$ filtered by subbundles. Since the essential image of $\nabla^{\mathrm{H}}_{\check{G}}$ are connections, they are smooth objects in $\rr{MTM}^{\rr{int}}_{\rr{resc}}(\bb{G}_{m,t})$. So, we can consider the composition
$$
\rr{Rep}_{\check{G}}\xrightarrow{\nabla^{\rr{H}}} \rr{MTM}^{\rr{int}}_{\rr{sm,resc}}(\bb{G}_{m,t}) \xrightarrow{\rr{Irr Fil}} \rr{Fil}_{\bb{G}_{m,t}}.
$$

Similar to \cite[\S 3.2.4]{Lovering17}, we call a functor $F\colon \rr{Rep}_{\check{G}}\to \rr{Fil}_{\bb{G}_m}$ an \textit{$\eta/k$-filtration} for some dominant cocharacter $\eta$ and some integer $k$, if for each representation $V$ of $\check G$, the decreasing filtration by subbundles $F(V)$ is given by
    $$F^\alpha \nabla_{\check{G}}(V):= \cc{O}_{\bb{G}_m}\otimes_{\dC}\biggl(\bigoplus\nolimits_{i \geq k\cdot \alpha} V_i\biggr),$$
where $V=\bigoplus_{i\in \dZ}V_i$ is the decomposition of $V$ induced by the $\bb{G}_m$-action $\eta_V:\bb{G}_m\xrightarrow{\eta} \check G\to \rr{GL}(V)$.

\begin{definition}\label{def:FiltGivenCocarac}
    We say that the irregular Hodge filtrations of $\nabla^{\rr{H}}_{\check{G}}$ are determined by $\eta/k$ for some dominant cocharacter $\eta\colon \bb{G}_{m}\to \check{G}$ and a positive integer $k$ if the composition $\rr{IrrFil}\circ\nabla^{\rr{H}}$ is an $\eta/k$-filtration. We say that the irregular Hodge numbers of $\nabla_{\check{G}}^{\rr{H}}$ are determined by $\eta/k$  if for any $V\in \rr{Rep}_{\check{G}}$ the ranks of $\rr{gr}^\alpha(\rr{IrrFil}\circ\nabla^{\rr{H}}(V))$ equal the
    dimensions of the eigenspaces $V_{k\cdot \alpha}$ of the operator $\eta(V)$.
\end{definition}
 If the irregular Hodge filtration of $\nabla^{\rr{H}}_{\check{G}}$ is determined by $\eta/k$, then so are its irregular Hodge numbers. The converse is not true in general.
Using these notions, we can summarize the above discussion as follows.
\begin{proposition}\label{lem:existence}
    The Frenkel--Gross connection $\nabla_{\check G}$ underlies a unique rescalable integrable mixed twistor $\scr{D}$-module $\nabla_{\check G}^{\rr{H}}$ of weight $0$, i.e.~a tensor functor
    \[
        \nabla_{\check{G}}^{\rr{H}}\colon \rr{Rep}_{\check{G}} \to \rr{MTM}_{\rr{resc}}^{\rr{int}}(\bb{G}_{m,t}).
    \]
    In particular, there is a uniquely defined irregular Hodge filtration on $\nabla_{\check{G}}$, and it is determined by $\eta/k$ for some cocharacter $\eta$ and some integer $k$.
\end{proposition}
\begin{proof}
    We already constructed rescalable integrable mixed twistor $\scr{D}$-module $\nabla_{\check{G}}^{\rr{H}}$ on $\nabla_{\check{G}}$ in \eqref{eq:EMHM-tensor}. To show the uniqueness of $\nabla_{\check{G}}^{\rr{H}}$, we decompose any representation $W$ of $\check{G}$ as $W=\bigoplus W_i\otimes M_i$ with $W_i$ irreducible and $M_i$ the multiplicities. By the functoriality  of $\nabla_{\check{G}}^{\rr{H}}$, we have
     \[
        \nabla_{\check G}^{\rr{H}}(W)=\bigoplus \nabla_{\check G}^{\rr{H}}(W_i)\otimes M_i^{\rr{H}},
     \]
     with $M_i^{\rr{H}}$ being trivial. As in \cref{lem:irreducible}, there is at most one rescalable integrable mixed twistor $\scr{D}$-modules of weight $0$ on an irreducible connection on $\bb{G}_m$. So each $\nabla_{\check G}^{\rr{H}}(W_i)$ is unique. Therefore, $\nabla_{\check G}^{\rr{H}}(W)$ is also unique.

    Lastly, since each $\nabla_{\check{G}}(V)$ is a connection on $\bb{G}_{m,t}$, its associated rescalable integrable mixed twistor $\scr{D}$-module gives rise to an irregular Hodge filtration (by subbundles) in a functorial way. In other words, we have a tensor functor from $\rr{Rep}_{\check{G}}$ to ($\frac{1}{k}\dZ$-indexed) filtered vector bundles on $\bb{G}_{m,t}$ for some integer $k$. By \cite[Lem. 3.3.1]{Lovering17}, this functor is an $\eta/k$-filtration for some cocharacter $\eta$.
\end{proof}

\subsection{Inverse Fourier transforms of  \texorpdfstring{$\widetilde{\nabla}_{\check{G}}(V)$}{}} \label{sec:inverse-fourier}

In this section, we study the Hodge properties of the (inverse) Fourier transform of the Frenkel-Gross connections and their variants. We first work in a slightly more general geometric framework to simplify our proof. The situation will be specified in the case of Frenkel--Gross connections at the end of this section.

\paragraph{The setup.}
Let $Y$ be a quasi-projective variety over $\dC$ and $g$ a regular function on $Y$. We consider the product $Y\times \bb{A}^1_{t}$ and we denote by $\rr{pr}_t$ resp.~ $\rr{pr}_Y$ the projection $Y\times \bb{A}^1_{t}\rightarrow \bb{A}^1_{t}$ resp.~the projection $Y\times \bb{A}^1_{t}\rightarrow Y$. We have the following diagram:
\[
\begin{tikzcd}[row sep=5pt]
    & Y\times \bb{A}^1_{t}\ar[ld,"t\cdot g"'] \ar[rd,"\rr{pr}_t"] &\\
    \bb{A}^1&  & \bb{A}^1_{t}.
\end{tikzcd}
\]
We start by collecting a few general facts about Fourier transforms of (complexes of) $\cD_{\dA^1_t}$-modules that we will use later.
Recall from the introduction that the Fourier transform of an object $\cM^\bullet \in D^b_h(\bb{A}^1_t)$ denoted by $\rr{FT}^{\pm 1} (\mathcal{M}^\bullet)$, is the complex of holonomic $\mathscr{D}_{\bb{A}^1_s}$-modules defined by
\[
   \mathrm{pr}_{s+} (\mathrm{pr}_t^+\mathcal{M}^\bullet\otimes \mathcal{E}^{\pm t\cdot s}),
\]
where $\mathrm{pr}_t$ and $\mathrm{pr}_{s}$  are the projections from $\bb{A}^1_t\times \bb{A}^1_s$ to $\bb{A}^1_t$ and $\bb{A}^1_s$ respectively. Similarly, we can define the Fourier transform of $\cM^\bullet$ with compact support by
\[
    \rr{FT}^{\pm}_{\dagger}(\cc{M}^\bullet):=\mathrm{pr}_{s\dagger} (\mathrm{pr}_t^+\mathcal{M}^\bullet\otimes \mathcal{E}^{\pm t\cdot s}).
\]
The two kinds of Fourier transforms are interchanged by the duality functor
\begin{equation}\label{eq:FTDual-1}
    \bb{D}\circ \rr{FT}^{-1}=\rr{FT}_\dagger\circ\bb{D},
\end{equation}
and the forget-support morphism induces an isomorphism $\rr{FT}_\dagger(\cc{M})\xrightarrow{\sim} \rr{FT}(\cc{M})$, see \cite[App.~2 Prop.~1.7]{malgrange91}, so that we also have
\begin{equation}\label{eq:FTDual-2}
\bb{D}\circ \rr{FT}^{-1}=\rr{FT} \circ\bb{D}.
\end{equation}
Moreover $\rr{FT}^{-1}\circ \rr{FT}=\rr{id}$.

\begin{proposition}\label{prop:inverse-fourier}
For any $N\in \textup{Mod}_h(\cD_Y)$ we have isomorphisms in $D^b(\cD_{\dA^1_s})$:
    \begin{equation*}
        \rr{FT}^{-1}(\rr{pr}_{t+}( \pr_Y^+N\otimes \cc{E}^{t\cdot g})) \cong g_+N
    \end{equation*}
    and
    \begin{equation*}
        \rr{FT}^{-1}(\rr{pr}_{t\dagger}( \pr_Y^+N  \otimes \cc{E}^{t\cdot g})) \cong g_\dagger N,
    \end{equation*}
here we see $g$ as a morphism $g\colon X \to \dA^1_s$.
\end{proposition}
\begin{proof}
    It suffices to calculate $\rr{FT}(g_+N)$ and $\rr{FT}(g_\dagger N)$. Consider the following diagram:
\[
\begin{tikzcd}
    & Y\times \bb{A}^1_t \arrow[ld, "\pr_Y"] \arrow[rr, "(g\times \rr{id})"] \arrow[rrrd, "\pr_t", bend left=50] &            & \bb{A}^1_s\times \bb{A}^1_t \arrow[ld, "\pr_s"] \arrow[rd, "\pr_t"'] &            \\
    Y \arrow[rr, "g"] & & \bb{A}^1_s & & \bb{A}^1_t
\end{tikzcd}\]
Then we have
\begin{align*}
    \mathrm{FT}\left(g_{+}N \right)
    &=\pr_{t+}\left(\pr_s^{+} g_{+} N \otimes  \mathcal{E}^{s \cdot t}\right) \\
    & =\pr_{t+}\left( \left(g\times  \mathrm{id}_{t}\right)_{+} \pr_{Y}^{+} N \otimes \mathcal{E}^{s \cdot t}\right)
    &&\text { smooth base change } \\
    & =\pr_{t+} \circ\left(g\times \mathrm{id}_{t}\right)_{+}\left(\pr_{Y}^{+} N \otimes \left(g\times  \mathrm{id}_{t}\right)^{+} \mathcal{E}^{s \cdot t}\right)
    &&\text { projection formula } \\
    & =\pr_{t+} \circ\left(g\times  \operatorname{id}_{t}\right)_{+}\left(\pr_{Y}^{+} N \otimes  \mathcal{E}^{g \cdot t}\right)
    &&\text { definition of } \mathcal{E}^{s \cdot t} \\
    & =\pr_{t+}\left(\pr_{Y}^{+} N\otimes \mathcal{E}^{g \cdot t}\right).
\end{align*}
In the same way, we obtain
\begin{equation}\label{eq:FL-direct-image}
    \mathrm{FT}^{-1}\left(g_{+}N \right)
    =\pr_{t+}\left(\pr_{Y}^{+} N\otimes \mathcal{E}^{-g \cdot t}\right).
\end{equation}

Using the usual commutation rules of the duality functor with direct and inverse images, as well as the commutation of duality with Fourier transformation discussed above, we derive from the previous computation the following:
\begin{align*}
\mathrm{FT}(g_\dag \bD_Y(N)) & =\mathrm{FT}(\bD_{\dA^1_s}(g_+ N)) \\
&= \bD_{\dA^1_t} (\mathrm{FT}^{-1}(g_+ N)) && \text{use formula \eqref{eq:FTDual-2}}\\
&=\bD_{\dA^1_t} \pr_{t+}\left(\pr_{Y}^{+} N\otimes \mathcal{E}^{-g \cdot t}\right)&& \text{use formula \eqref{eq:FL-direct-image}}\\
&= \pr_{t\dag} \bD_{Y\times\dA^1_t}\left(\pr_{Y}^{+} N\otimes \mathcal{E}^{-g \cdot t}\right) \\
&= \pr_{t\dag} \left(\pr_{Y}^{+} \bD_Y N\otimes \bD_{Y\times\dA^1_t}\mathcal{E}^{-g \cdot t}\right)\\
&= \pr_{t\dag} \left(\pr_{Y}^{+} \bD_Y N\otimes \mathcal{E}^{g \cdot t}\right) &&\text{$\cE^{-g\cdot t}$ is locally free.}\\
\end{align*}
Since this holds for any $N$, we obtain $\mathrm{FT}(g_\dag N) =
\pr_{t\dag} \left(\pr_{Y}^{+} N\otimes \mathcal{E}^{g \cdot t}\right)$, as required.
\end{proof}

Let $j\colon \bb{G}_{m,t}\hookrightarrow \bb{A}^1_{t}$ be the inclusion and let $N$ be a $\cD_Y$-module. Consider
\[
\cF:= j^+\pr_{t+}(\cc{E}^{t\cdot g}\otimes \pr_Y^+ N),
\]
which is a priori an object in $D^b(\scr{D}_{\bb{G}_{m,t}})$.
Recall that for $a\in\dQ$ and $\chi=\exp(2\pi i a)$, we write $\bb{L}_{\chi}:=(\mathcal{O}_{\bb{G}_{m,t}},\mathrm{d}+a\frac{\mathrm{d}t}{t})$ and we put $\cF_\chi:=\cF\otimes \bb{L}_\chi$.
We will consider the following geometric situation:
\begin{assumption}\label{assum:OneToFour}
We take $a\in\dQ\setminus\dZ$, i.e.~$\chi\neq1$, and subject the above objects to the following four conditions.
\begin{enumerate}
    \item $\cF$ is a $\scr{D}_{\G_{m,t}}$-module concentrated in degree $0$, regular singular at $0$, and irregular singular at $\infty$ with slopes $1$. Notice that the same assumption then automatically holds for $\cF_\chi$.

    \item The forget support morphism
    \begin{equation}\label{eq:assumption-1}
        j^+\mathrm{pr}_{t\dagger}(\mathcal{E}^{\pm t\cdot g} \otimes \pr_Y^+ N) \to j^+\mathrm{pr}_{t+}(\mathcal{E}^{\pm t\cdot g} \otimes \pr_Y^+ N)=\cF
    \end{equation}
    induces an isomorphism of $\scr{D}_{{\bb{G}_{m,t}}}$-modules concentrated in degree $0$.

    \item The morphism
    \begin{equation}\label{eq:assumption-2}
        j_\dagger (\cF_\chi)
        \to
        j_+(\cF_\chi)
    \end{equation}
    is an isomorphism, so in particular both $j_\dagger (\cF_\chi)$ and $ j_+(\cF_\chi)$ are isomorphic to $j_{\dagger+}(\cF_\chi)$.

    \item $N$ underlies a pure Hodge module $N^{\rr{H}}\in \rr{MHM}(Y)$ of weight $\dim Y$.
\end{enumerate}

\end{assumption}

We continue our discussion of various Fourier transforms. We restrict, however, from now on our attention to objects satisfying the above four conditions.

For two $\scr{D}_{\bb{A}^1}$-modules $M$ and $N$, the additive convolutions are $\scr{D}_{\bb{A}^1}$-modules defined as
\[
    M\star_+N:=\rr{sum}_+(M\boxtimes N)
    \quad \text{and} \quad
    M\star_\dagger N:= \bb{D}(\bb{D}(M)\star_+\bb{D}(N)),
\]
where $\rr{sum}\colon \bb{A}^1\times \bb{A}^1\to \bb{A}^1$ is the summation map.

\begin{proposition}\label{prop:inverse-Fourier-twisted}
    Under \cref{assum:OneToFour} we have
    \begin{equation}\label{eq:inverse-fourier}
        \rr{FT}^{-1}(j_+\cF_\chi)=g_+N\star_+j_+\bb{L}_{\chi^{-1}}
    \end{equation}
    and
    \begin{equation}\label{eq:inverse-fourier-2}
        \rr{FT}^{-1}(j_\dagger\cF_\chi)=g_\dagger N\star_\dagger j_\dagger\,\bb{L}_{\chi^{-1}}.
    \end{equation}
    Here $\star_+$ and $\star_\dagger$ are the additive convolutions.
\end{proposition}

\begin{proof}
Since the Fourier transform interchanges the tensor product and the additive convolution, by \cite[Equation (1.1.2)]{DS13} and the first equation in \cref{prop:inverse-fourier}, we have
\[
    \begin{split}
        \rr{FT}^{-1}&(\mathrm{pr}_{t+}(\pr_{Y}^+N\otimes \mathcal{E}^{t\cdot g})\otimes j_+\bb{L}_\chi)\\
        &=\rr{FT}^{-1}(\mathrm{pr}_{t+}(\pr_{Y}^+N\otimes\mathcal{E}^{t\cdot g}))\star_+\rr{FT}^{-1}( j_+\bb{L}_\chi)\\
        &=g_+N\star_+j_+\bb{L}_{\chi^{-1}}.
    \end{split}
\]
Notice that the adjunction morphism $\rr{id}\to j_+j^+$ induces a morphism
\[
    \mathrm{pr}_{t+}(\pr_{Y}^+N\otimes\mathcal{E}^{t\cdot g})\to j_+\cF
\]
whose kernel and cokernel are supported at $0$. Hence, we deduce an isomorphism
\[
    \mathrm{pr}_{t+}(\pr_{Y}^+N\otimes\mathcal{E}^{t\cdot g})\otimes j_+\bb{L}_\chi\simeq j_+\cF_\chi.
\]
Combining the above identities, we deduce \eqref{eq:inverse-fourier}.

Similarly, we have
\[
    \rr{FT}(\mathrm{pr}_{t+}(\pr_{Y}^+N\otimes\mathcal{E}^{-t\cdot g})\otimes j_+\bb{L}_{\chi^{-1}})
    =\rr{FT}(\mathrm{pr}_{t+}(\pr_{Y}^+N\otimes\mathcal{E}^{-t\cdot g}))\star_+\rr{FT}( j_+\bb{L}_{\chi}^{-1})
    =g_+N\star_+j_+\bb{L}_{\chi}.
\]
Notice that
\[
    \rr{pr}_{t+}(\pr_{Y}^+N\otimes\cc{E}^{-t\cdot g})|_{\bb{G}_{m}}\simeq
    \rr{pr}_{t\dagger }(\pr_{Y}^+N\otimes \cc{E}^{-t\cdot g})|_{\bb{G}_{m}}=\bb{D}_{\bb{G}_{m}}(\cF)=\cF^\vee,
\]
we argue similarly above that
\[
    \rr{pr}_{t+}(\pr_{Y}^+N\otimes\cc{E}^{-t\cdot g})\otimes j_+\bb{L}_{\chi^{-1}}=j_+(\cF_{\chi})^\vee.
\]
So, we deduce that
\[
    \rr{FT}(j_+(\cF_{\chi })^\vee)=g_+N\star_+j_+\bb{L}_{\chi}
\]

Applying the duality functor, we have
\[
    \begin{split}
        \rr{FT}^{-1}_\dagger(j_\dagger \cF_{\chi})
        &=\bb{D}_{\bb{A}^1}\circ \rr{FT} (j_+(\cF_{\chi})^\vee)
                =\bb{D}_{\bb{A}^1}(g_+ N\star_+ j_+\bb{L}_{\chi})\\
        &=\bb{D}_{\bb{A}^1}(\bb{D}_{\bb{A}^1}(g_\dagger N)\star_+ \bb{D}_{\bb{A}^1}(j_\dagger\bb{L}_{\chi^{-1}}))
                =g_\dagger N\star_\dagger j_\dagger\bb{L}_{\chi^{-1}},
    \end{split}
\]
which is \eqref{eq:inverse-fourier-2}.
\end{proof}
\paragraph{Mixed Hodge modules on the inverse Fourier transformation}
Notice that $ g_\dagger N\star_\dagger j_\dagger\bb{L}_{\chi}$ resp.~$g_+ N\star_+j_+\bb{L}_{\chi}$ underlies a mixed Hodge module
\[
    \rr{sum}_{!}\bigl( g_! N^{\rr{H}}\boxtimes j_!\bb{L}_\chi^{\rr{H}}\bigr)
    \quad
    \text{resp.}
    \quad \rr{sum}_{*}\bigl(g_* N^{\rr{H}}\boxtimes j_*\bb{L}_\chi^{\rr{H}}\bigr),
\]
where $N^{\rr{H}}$ is the pure Hodge module on $Y$ of weight $\dim Y$ with underlying $\scr{D}$-module $N$, and $\bb{L}_\chi^{\rr{H}}$ is the pure complex Hodge module on $\bb{G}_{m}$ of weight $1$ with underlying $\scr{D}$-module $\bb{L}_\chi$.
\begin{lemma}\label{lem:ConvPureHodge} Under \cref{assum:OneToFour}, the $\cD_{\dA^1}$-module $g_+N\star_+j_+\bb{L}_{\chi}$ underlies a pure Hodge module of weight $\dim( Y)+1$ on $\bb{A}^1$.
\end{lemma}
\begin{proof}
It follows from \eqref{eq:assumption-1}, \eqref{eq:assumption-2}, \eqref{eq:inverse-fourier}, and \eqref{eq:inverse-fourier-2} that the forget-support morphism
\begin{equation}\label{eq:forget-support}
    g_\dagger N\star_\dagger j_\dagger\bb{L}_{\chi} \to g_+N\star_+j_+\bb{L}_{\chi},
\end{equation}
being the composition
\[
    \rr{FT}^{-1}(j_\dagger \cF_\chi)\to
    \rr{FT}^{-1}(j_{\dagger+} \cF_\chi)  \to
    \rr{FT}^{-1}(j_{ +} \cF_\chi),
\]
is an isomorphism.

By the properties of conservation of weights of mixed Hodge modules (see \cite[Prop. 1.7]{Saito-intro}), the two mixed Hodge modules $\rr{sum}_{!}\bigl( g_! N^{\rr{H}}\boxtimes j_!\bb{L}_\chi^{\rr{H}}\bigr) $ and $\rr{sum}_{*}\bigl(g_* N^{\rr{H}}\boxtimes j_*\bb{L}_\chi^{\rr{H}}\bigr)$ are mixed of weight $\leq \dim Y+1$ and $\geq \dim Y+1$ respectively. Since \eqref{eq:forget-support} can be lifted to an isomorphism of mixed Hodge modules, we conclude that they are pure Hodge modules of weight $\dim(Y)+1$.
\end{proof}

\paragraph{The case of Frenkel--Gross connections}
Let $\check{G}$ be a simply connected simple group (i.e., $G$ an almost simple group of adjoint type) and $V$ a representation of $\check{G}$. We use notation from \eqref{eq:HNY-diagram} and consider the case that $X=\rr{Gr}^\circ$, $g=F(x,1)$, $\pr_t$ is the projection to $\bb{A}^1_t$ such that $\rr{pr}_t\circ (\rr{id}\times j)=j\circ \pi^\circ$, $N=\rr{IC}_{V}|_{\rr{Gr}^\circ}$ as follows:
\begin{equation}
\begin{tikzcd}
     &\rr{Gr}^\circ \times \bb{G}_{m} \ar[rd,"\pi^\circ"] \ar[ld,"t\cdot g"'] \ar[r,hook,"\rr{id}\times j"]& \rr{Gr}^\circ \times \bb{A}^1_{t} \ar[rd,"\rr{pr}_t"] & \\
        \bb{A}^1 &&\bb{G}_{m} \ar[r,hook, "j"]& \bb{A}^1_{t} .
\end{tikzcd}.
\end{equation}
From the discussion above, we deduce the following corollary:
\begin{corollary}\label{lem:inverse-fourier}
    The inverse Fourier transform $\mathrm{FT}^{-1}\circ j_+( \widetilde{\nabla}_{\check{G}}(V)\otimes \bb{L}_\chi)$ is a regular holonomic $\scr{D}_{\bb{A}^1}$-module. If $\chi$ is non-trivial, the inverse Fourier transform can be lifted to a pure Hodge module $\widetilde M^{\mathrm{H}}_\chi$ on ${\bb{A}^1}$ of weight $w:=\dim(\rr{Gr}_{\lambda})+1=\langle 2\rho,\lambda\rangle+1$.
\end{corollary}
\begin{proof}
First note that $t\cdot g$ can be extended to $\mathrm{Gr}^\circ \times \bb{A}^1_t$. Recall that  $\rr{IC}_{V,\rr{GR}^\circ}=\rr{pr}_Y^+\rr{IC}_{V,\rr{Gr}^\circ}$ via the projection $\rr{pr}_Y\colon \rr{GR}^\circ=\rr{Gr}^\circ \times \bb{G}_m\to \rr{Gr}^\circ$.
It then follows from \cref{prop:HNY-D-module} that
\[
    \widetilde{\nabla}_{\check{G}}(V)={\pi}^\circ_+(\cc{E}^{t\cdot g}\otimes \rr{pr}_Y^+\rr{IC}_{V,\rr{Gr}^\circ}).
\]
It can be verified that the conditions of \cref{assum:OneToFour}, as outlined at the beginning of \cref{sec:inverse-fourier}, are satisfied for this choice of $(Y,g,N,\cF= \widetilde{\nabla}_{\check{G}}(V)\otimes \bb{L}_\chi)$. Indeed, given our choices, the first and last conditions are straightforward. The second condition is verified by \eqref{eq:!vs*}, and the third condition holds because $\nabla_{G}$ is regular singular and has principal unipotent monodromy at $0$.

It then follows from \cref{prop:inverse-Fourier-twisted} and \cref{lem:ConvPureHodge} that $\mathrm{FT}^{-1}\circ j_+( \widetilde{\nabla}_{\check{G}}(V)\otimes \bb{L}_\chi)$ underlies a pure Hodge module on ${\bb{A}^1}$ of weight $\langle 2\rho,\lambda\rangle+1$.
\end{proof}

\begin{corollary}\label{cor:unicity}
Assume that $\nabla_{\check{G}}(V)$ is an irreducible connection on $\G_{m,t}$. Then there is a unique rescalable integrable mixed twistor $\scr{D}$-modules of weight $0$ on $\dA^1_t$ with underlying $\cD_{\dA^1_t}$-module equal to $j_{\dag+}\tilde{\nabla}_{\check G}(V)$.
Furthermore, $\mathbf{FT}\circ v\bigl(\widetilde M^{\rr{H}}_\chi\bigr) $ coincides with the rescalable mixed twistor $\cD$-module up to a Tate twist defined in \eqref{eq:EMHM-chi}, where $\mathbf{FT}$ is the Fourier transform for $\mathrm{MTM}^{\rr{int}}(\bb{A}^1)$.
\end{corollary}
\begin{proof}
Apply \cref{lem:irreducible} for $Y=\dA^1$ and for $N=j_{\dag+}\tilde{\nabla}_{\check G}(V)$ as well as for $N=j_{+}(\tilde{\nabla}_{\check G}(V)\otimes \dL_\chi)$.
\end{proof}

\subsection{The irregular Hodge numbers of \texorpdfstring{$\nabla_{\check{G}}$}{}} \label{sec:proof}

Now we prove \cref{thm:main}.

\begin{lemma}\label{lem:reduction1}
    Let $G_1$ and $G_2$ be two isogenous almost simple groups. If the irregular Hodge filtrations of $\nabla_{\check G_1}^{\rr{H}}$ are determined by $\rho_{G_1}$, then the irregular Hodge filtrations of $\nabla_{\check G_2}^{\rr{H}}$ are determined by $\rho_{G_2}$.
\end{lemma}
\begin{proof}
    Without loss of generality, we may assume that we have an isogeny $\pi\colon \check{G}_1\to \check{G}_2$ or $\pi\colon \check{G}_2\to \check{G}_1$. In the former case, we deduce the irregular Hodge filtrations of $\nabla_{\check{G}_2}^{\rr{H}}$ by the ``trivial'' functoriality, as discussed in \cref{exe:functoriality}.

    In the latter case, assume that the irregular Hodge filtration of $\nabla_{\check{G}_2}^{\rr{H}}$ is determined by $\eta/k$. By the ``trivial'' functoriality again, $\eta/k$ will induce $2\rho_{G_1}$ on $\check{G}_1$, i.e.~$\pi\circ \eta^2=(2\rho_{G_1})^k$.
    Hence, the cocharacter $\pi\circ( \eta^2\cdot (2\rho_{G_2})^{-k})$    is trivial, which forces $\eta^2\cdot (2\rho_{G_2})^{-k}$ to take values in $\ker \pi$. Since there is no nontrivial homomorphism from $\bb{G}_{m}$ to finite groups, we conclude that $\eta^2= (2\rho_{G_2})^{k}$, which means that $\eta/k=\rho_{G_2}$.
\end{proof}

\begin{lemma}\label{lem:reduction2}
    If the irregular Hodge numbers of $\nabla_{\check G}^{\rr{H}}(V)$ are given by $\rho(V)$ for one faithful representation $V$ of ${\check G}$, then the irregular Hodge numbers of $\nabla_{\check G}^{\rr{H}}$ are given by $\rho$.
\end{lemma}
\begin{proof}
As algebraic vector bundles on $\bb{G}_m$ are trivial, we can (non-canonically) identify any basis of $\nabla_{\check{G}}(L)$ with a basis of $L$ for any representation $L$. In particular, the irregular Hodge filtration on $\nabla_{\check{G}}(L)$ determines a flag (indexed by $\frac{1}{k}\dZ$)
\[
    \scr{F}_{\alpha/k}\subset \ldots\subset \nabla_{\check{G}}(L)
\]
for $\alpha\in\dZ$, such that $h^\alpha=\rk \scr{F}_{\alpha/k}/\scr{F}_{(\alpha-1)/k}$.

By \cref{lem:existence}, the irregular Hodge filtration is determined by $\eta/k$ for some cocharacter $\eta$ of $\check{G}$ and an integer $k$. It follows that $\eta(V)$ and $\rho^k(V)$ determines the same flag of $V$. As $V$ is faithful, $\eta$ and $\rho^k$ correspond to the same parabolic subgroup $\check{P}$ of $\check{G}$ up to conjugacy. It follows that for any representation $W$ of $\check{G}$, $\eta(W)$ and $\rho^k(W)$ determine again the same flag on $W$. Hence, the irregular Hodge numbers of $\nabla_{\check{G}}(W)$ are determined by $\rho(W)$.
\end{proof}

Now, we need to calculate the irregular Hodge numbers of some specific representations $V$. Similar to \cite[\S5.2.1]{Qin23}, we consider the case where the local monodromy of $\nabla_{\check{G}}(V)$ at $0$ consists of Jordan blocks of different sizes. This assumption ensures that the primitive parts of the Lefschetz decomposition of some limiting Hodge structures are of rank $1$, hence of Hodge--Tate type. If the primitive parts have dimensions bigger than $1$, there are more possibilities. In this case the computation of the irregular Hodge numbers may be more involved. As we will see in \cref{lem:different-Jordan-blocks} below, this case will not occur in our investigation.

\begin{proposition}\label{prop:Hodge-numbers}
    Assume that $V$ is irreducible and the local monodromy of $\nabla_{\check{G}}(V)$ at $0$ consists of Jordan blocks of sizes $r_1< r_2< \dots< r_k$. Then the nearby cycle
 $\psi_{1/\tau}\widetilde{M}_{\chi}^{\rr{H}}$ is of Hodge--Tate type, and the Hodge numbers coincide with the irregular Hodge numbers of $ \nabla_{\check{G}}(V)$ up to a global shift.

 More precisely, we have for $p\in\frac{1}{2}\dZ$ that
    \[
        h^p_{\rr{irr}}=\#\{(i,a)\mid 2p=r_i-2a,~0\leq a\leq r_i \}.
    \]
\end{proposition}
\begin{proof}
    Let $\widetilde{\nabla}_{\check{G}}(V)=[h]^+\nabla_{\check{G}}(V)$, $\bb{L}_{\chi}$ and $g$ be as in the previous section. In order to prove the proposition, it suffices to calculate the irregular Hodge numbers of $\widetilde{\nabla}_{\check{G}}(V)\otimes \bb{L}_\chi$ for any $\chi\neq 1$.

Recall that $j\colon \bb{G}_m\hookrightarrow \bb{A}^1$ denotes the canonical embedding. Notice that $\widetilde\nabla_{\check{G}}(V)\otimes\bb{L}_{\chi}$ is regular at $0$, with quasiunipotent monodromy of eigenvalue $\chi$. So we have
\[
    j_{\dagger}( \widetilde{\nabla}_{\check{G}}(V)\otimes\bb{L}_{\chi})=j_{\dagger+}( \widetilde{\nabla}_{\check{G}}(V)\otimes\bb{L}_{\chi})=j_{+}( \widetilde{\nabla}_{\check{G}}(V)\otimes\bb{L}_{\chi}).
\]

Recall also that $\widetilde M_{\chi}:=\rr{FT}^{-1}j_{\dagger+}( \widetilde{\nabla}_{\check{G}}(V)\otimes\bb{L}_{\chi})$ is the inverse Fourier transform of $j_{\dagger+}( \widetilde{\nabla}_{\check{G}}(V)\otimes\bb{L}_{\chi})$, which can be lifted to a pure Hodge module $\widetilde M^{\mathrm{H}}_\chi$ by \cref{lem:inverse-fourier}.

Let $\psi_t\,( \widetilde{\nabla}_{\check{G}}(V)\otimes \bb{L}_\chi)$ be the nearby cycle module at $0$ of $j_{+}( \widetilde{\nabla}_{\check{G}}(V)\otimes \bb{L}_\chi)$. Notice that the local monodromy of $ \widetilde{\nabla}_{\check{G}}(V)$ at $0$ is also unipotent with Jordan blocks of sizes $r_1<\ldots <r_k$. So that of $\psi_t\,( \widetilde{\nabla}_{\check{G}}(V)\otimes \bb{L}_\chi)$ is quasiunipotent with eigenvalue $\chi$ and with Jordan blocks of sizes $r_1<\ldots<r_k$.

By construction, we have $\widetilde{M}_\chi=\mathrm{FT}^{-1} j_{\dagger+}( \widetilde{\nabla}_{\check{G}}(V)\otimes \bb{L}_\chi)$. Applying the (inverse) stationary phase formula (cf.~\cite[Prop.~VII.2.4]{malgrange91} or \cite[Prop.~2.3]{Sab06hyp}), it follows that
\begin{equation}\label{eq:stationary-phase}
        \psi_t\,( \widetilde{\nabla}_{\check{G}}(V)\otimes \bb{L}_\chi)=\phi_{1/\tau}\widetilde M_\chi=\psi_{1/\tau,-1}\widetilde M_\chi
\end{equation}
and the corresponding nilpotent part of the monodromy operator $N$ has Jordan blocks of sizes $r_1<\dots<r_k$ respectively. In particular, $\rk\,\widetilde M_\chi=r$.

By our assumption, the (non-zero) primitive parts of the Lefschetz decomposition of $\psi_{1/\tau}\widetilde M^{\rr{H}}_\chi$, denoted by $P_{r_1},\ldots,P_{r_k}$, are one-dimensional and of Hodge--Tate type. Since $\widetilde M^{\rr{H}}_\chi$ is pure of weight $w=\langle 2\rho,\lambda \rangle+1$, the primitive parts $P_{r_1},\dots,P_{r_k}$ are pure of weights $r_1+w-1,\ldots,r_k+w-1$ respectively. In other words,
\[
    P_{r_i}=\mathrm{gr}^W_{r_i+w-1}P_{r_i}
\]
for $1\leq i\leq k$. Then, by the Lefschetz decomposition, for each $\ell\in\dZ$, the graded quotient $\mathrm{gr}^W_{\ell}\psi_{1/\tau}\widetilde M^{\mathrm{H}}_\chi$ is Hodge--Tate
of dimension
\[
    \#\{(i,a)\mid \ell=r_i+w-1-2a, 0\leq a\leq r_i \}.
\]
Moreover, we have
	\begin{equation*}
		\begin{split}
			\dim \rr{gr}_F^p \psi_{1/\tau}\widetilde M^{\mathrm{H}}_\chi=\dim \mathrm{gr}^W_{2p}\psi_{1/\tau}\widetilde M^{\rr{H}}_\chi.
		\end{split}
	\end{equation*}

Since the eigenvalues of $\widetilde{M}_\chi^{\rr H}$ at $\infty$ are different from $1$, we conclude using \cite[(7)]{sabbah-hypergeometric} that
\begin{center}
        $\dim \mathrm{gr}_{F_{irr}}^{p+a}( \widetilde{\nabla}_{\check{G}}(V)\otimes \bb{L}_\chi)=\dim \mathrm{gr}_{F}^{p}\psi_{1/\tau,-1} \mathbf{FT}^{-1}(\widetilde{\nabla}_{\check{G}}(V)\otimes \bb{L}_\chi)$
\end{center}
which by \cref{cor:unicity} coincide with
\[
    \dim \mathrm{gr}^p_F\psi_{1/\tau,-1}\widetilde{M}_{\chi}^{\rr{H}}=\#\{(i,a)\mid 2p=r_i+w-1-2a, 0\leq a\leq r_i \}
\]
up to a global shift, where $a$ is a rational number such that $\chi=\exp(2\pi i a)$. Notice that the irregular Hodge numbers of $ \widetilde{\nabla}_{\check{G}}(V)\otimes \bb{L}_\chi$ agree with those of $ \widetilde{\nabla}_{\check{G}}(V)$ (as well as $\nabla_{\check{G}}(V)$) up to a global shift. Therefore, we have shown the formula for the irregular Hodge numbers of $\nabla_{\check{G}}(V)$ as announced in the statement.
\end{proof}

\begin{corollary}\label{thm:minuscule}
    Assume that $\nabla_{\check{G}}(V)$ is irreducible and its local monodromy at $0$ consists of Jordan blocks of different sizes. Then the Hodge numbers of $\nabla_{\check{G}}(V)$ are determined by $\rho(V)$.
\end{corollary}
\begin{proof}
    Let $V=\bigoplus_{d\in  \frac{1}{2}\dZ} V_d$ be the degree decomposition of $V$ with respect to $2\rho\colon \bb{G}_{m}\to \check{G}$, such that $2\rho(t)$ acts on $V_d$ by multiplication by $t^{2d}$. There is a decreasing filtration induced by $\rho$, i.e.~the filtration defined by
    \[
        F_\rho^p:=\bigoplus_{p\leq 2d }V_d
    \]
    for $p\in\dZ$. Equivalently, $2\rho$ can be seen as an element in $\check{\mathfrak{g}}$, which induces an endomorphism $2\rho(V)$ on $V$ such that $v\in V_d$ if and only if $2\rho(V)\cdot v = 2dv$.

    Notice that one can upgrade the principal nilpotent operator $N$ to a $\mathfrak{sl}_2$-triple $(x,y,h)$ such that $y=N$. Recall that $\rho $ can equivalently be written as $\rho =\sum_{i=1}^n\omega_i$, the sum of \emph{fundamental coweights} of $\check{G}$ satisfying $\langle \al_i^\vee,\omega_j \rangle = \delta_{ij}$ for $\{\al_1^\vee,\ldots,\al_n^\vee\}=\Delta$ an ordering of the simple roots of $\check{G}$. Hence, we have
\begin{equation}\label{eq:[N,rho]=-N}
    [N,\rho]=\T\sum_{i,j=1}^n\bigl[X_{-\al_i^\vee},\omega_j] = \sum_{i,j=1}^n\langle-\al_i^\vee,\omega_j \rangle\,X_{-\al_i^\vee} = -N.
\end{equation}
So, we deduce that
\begin{equation}\label{eq:HTwoRho}
    h=2\rho.
\end{equation} In particular, we have
    \[
        \rho(V) N(V)v=N(V)\rho(V)v+N(V)v=(d+1)N(V)v
    \]
    for $v\in V_d$. Hence, one has for $d\in \frac{1}{2}\dZ$ that
    \[
        N(V)F_\rho^{2d}V\subset F_\rho^{2d+2}V.
    \]
    On the other hand, viewing $V$ as the representation of $\mathfrak{sl}_2$ via the composition $\mathfrak{sl}_2\to \check{\mathfrak{g}}\to \mathfrak{gl}(V)$, we have
    \[
        N(V)^{2d}\colon \rr{gr}_{F_\rho}^{-2d}V= V_{-d}\xrightarrow{\sim} V_d=\rr{gr}_{F_\rho}^{2d}V
    \]
    when $d\leq 0$ by the representation theory of $\mathfrak{sl}_2$ \cite[Thm.~7.2]{Humphreys72}. So $F_\rho^{2\bullet}$ is the same as the monodromy-weight filtration on $V$ with respect to $N(V)$.

    Notice that $\psi_t\,( \widetilde{\nabla}_{\check{G}}(V)\otimes \bb{L}_\chi)$ with its monodromy-weight filtration is abstractly isomorphic to $V$ with its monodromy-weight filtration with respect to $N(V)$,
     and so is $\psi_{1/\tau,-1}\widetilde{M}_{\chi}^{\rr{H}}$ by the identification from \eqref{eq:stationary-phase}. By \cref{prop:Hodge-numbers}, $\psi_{1/\tau,-1}\widetilde{M}_{\chi}^{\rr{H}}$ is of Hodge--Tate type and have Hodge numbers equal to the irregular Hodge numbers of $\nabla_{\check{G}}(V)$ up to a global shift, which are therefore determined by $\rho(V)$.
\end{proof}

To achieve the proof of \cref{thm:main}, it suffices to show the existence of representations $V$ such that the Jordan blocks of $N(V)$ have different sizes. In fact, adjoint representations are what we are looking for, as suggested by the following lemma:

\begin{lemma}\label{lem:different-Jordan-blocks}
    For a simple group $\check{G}$ not of type $D_{2n}$, the nilpotent operator $N$ acting on the adjoint representation $V=\mathfrak{g}$ has Jordan blocks of different sizes.
\end{lemma}
\begin{proof}
Let $V$ be the adjoint representation $V=\check{\mathfrak{g}}=\rr{Lie}(\check{G})$ of $\check{G}$ and $l$ be the rank of $\check{\mathfrak{g}}$. Notice that $N\in\check{\mathfrak{g}}$ is principally nilpotent, which can be enhanced into a principal $\mathfrak{sl}_2$. By \cite[\S6.5 \& Cor.~8.7]{Konstantprincipalthreedimensional59}, we can decompose $V=\check{\mathfrak{g}}$ into $\bigoplus_{i} V_{k_i}$, where $\dim V_{k_i}$ are $(2k_i+1)$-dimensional irreducible $\mathfrak{sl}_2$-representations, and the numbers $n_1\le\ldots\le n_l$ are the \textit{exponents} of $\check{\mathfrak{g}}$.

By \cite[Plate I-VII]{BourbakiBourbakiLiegroups02}, the exponents of $\check{\mathfrak{g}}$ are distinct numbers. So the dimensions of irreducible representations (as well as the Jordan blocks of $N$) are distinct numbers.
\end{proof}

\begin{proof}[{Proof of \cref{thm:main}}]
    The rescalable integrable mixed twistor $\scr{D}$-module on $\nabla_{\check{G}}$ exists by \cref{lem:existence}. To calculate its irregular Hodge numbers, it suffices to prove the theorem for one type in each row of diagram \eqref{eq:monodromy-groups} by the functoriality \eqref{eq:pushout}. So we assume that $\check{G}$ is of type $A_n$, $B_n$, $E_7$, $E_8$, and $F_4$.

    For each type as above, by \cref{lem:reduction1}, it suffices to prove the theorem for one almost simple group. So we prove the theorem for simple groups for each type.

    By \cref{lem:reduction2}, it suffices to know the Hodge numbers of $\nabla_{\check{G}}(V)$ for one faithful representation of $\check{G}$. So we choose $V$ as the adjoint representation of these groups.
    As $\check{G}$ is not of type $D_{2n}$, from \cref{lem:different-Jordan-blocks}, we deduce that the local monodromy at $0$ of the irreducible connection $\nabla_{\check{G}}(V)$ has Jordan blocks of different sizes. Therefore, we conclude the proof by \cref{thm:minuscule}.
\end{proof}

\section{On a conjecture of Katzarkov--Kontsevich--Pantev}\label{sec:KKP-conj}

As explained in the introduction, we are going to apply our main result \cref{thm:main} to the study of a conjecture by Katzarkov-Kontsevich-Pantev about irregular Hodge numbers of certain Landau-Ginzburg models for smooth projective Fano varieties.

By a \textit{Landau--Ginzburg model}, we mean a pair $(Y,w)$ consisting of a quasi-projective variety $Y$ and a regular function $w\colon Y\to \bb{A}^1$. In \cite[Conj.~3.7]{KKP2}, the authors defined \textit{Landau--Ginzburg Hodge numbers} of $(Y,w)$ as
\[
    f^{p,q}(Y,w):=\dim_\C \rr{H}^p(\mathsf{Z},\Omega^q_{\mathsf{Z}}(\log D,w)),
\]
where $Z$ is a smooth proper compactification of $Y$ such that $D:=Z\backslash Y$ is a simple normal crossing divisor, where $w$ extends to a morphism $\tilde{w}\colon Z \to\bb{P}^1$, and where $\Omega_Z^q(\log D,w)$ is the sheaf
\[
    \ker(\Omega_Z^q(\log D )\xrightarrow{\rr{d}w\wedge} \Omega_Z^{q+1}(\log D )).
\]
It is known \cite[Theorem 1.3.2]{Esnault-Sabbah-Yu} that we have
\[
f^{p,q}(Y,w)=\dim \gr^p_{F_{\textup{Yu}}} \rr{H}^{p+q}_{\rr{dR}}(\Omega_Y^\bullet,d+dw),
\]
where $F_{\textup{Yu}}^\bullet$ is the Yu filtration on the twisted de Rham cohomologies as defined in \cite{Yu14}. In particular, as the Yu filtration is independent of the choice of a compactification $Z$, so are the numbers $f^{p,q}(Y,w)$.

The authors of \cite{KKP2} predicted that when a Landau--Ginzburg model $(Y,w)$ is the mirror of a  Fano variety $X$, the Hodge numbers of $X$ are related to the Landau--Ginzburg Hodge numbers of $(Y,w)$ by the formula
\begin{equation}\label{eq:conj-KKP}
    f^{p,q}(Y,w) = h^{p,n-q}(X),
\end{equation}
where $n=\dim X$.

This conjecture has been verified in a number of cases, including the case of convenient and non-degenerate Laurent polynomials (\cite[Thm.~3.6]{Sabbah-KKP-conj}, where the variant \cite[Conj.~3.6]{KKP2} is shown).

We prove this conjecture when $X=\check{G}/\check{P}$ is a minuscule homogeneous space. The mirror Landau--Ginzburg model is $(Y=G\overset{\circ}{/}P,w)$, where $G$ resp.~$P$ are Langlands dual to $\check{G}$ resp.~$\check{P}$ (subject to the choice of a root datum), $G\overset{\circ}{/}P\subset G/P$ is the open projected Richardson variety, and $w$ is induced by some decoration function, see \cite[\S1.4]{Lam-Templier}.

\begin{theorem}\label{thm:Kontsevich}
    For $X=\check{G}/\check{P}$ a minuscule homogeneous space with $\dim(X)=n$, consider the Landau--Ginzburg model $(Y,w)$ from above.
    Then \eqref{eq:conj-KKP} holds, namely:
    \[
    f^{p,q}(Y,w)=h^{p,n-q}(X) = \begin{cases}
        0, &\text{if $q\neq n-p$,} \\
        h_{irr}^{p-\frac n2}, &\text{if $q=n-p$.}
    \end{cases}
    \]
    In particular, both numbers are determined by $\rho=2\rho/2$ up to a shift.
\end{theorem}
Before giving the proof, we explain how the Hodge numbers $f^{p,q}(Y,w)$ are related to our irregular Hodge numbers $h_{\rr{irr}}$ of Frenkel--Gross connections.
\begin{lemma}\label{lem:Yu-vs-irr}
    The numbers $f^{p,q}(Y,w)$ are zero if $q\neq n-p$ and are equal to $h_{\rr{irr}}^{p-n/2}$ when $q=n-p$, which is determined by $\rho$ up to a twist.
\end{lemma}
\begin{proof}
    Following \cite[\S1.3]{Lam-Templier}, let $X_{(G,P)}$ be the \textit{(parabolic) geometric crystal}, $f\colon X\to \bb{A}^1$ the \textit{decoration function}, and $\pi\colon X_{(G,P)}\to \bb{G}_m$ the \textit{highest weight function}. The fiber $X_t=\pi^{-1}(t)$ for $t\in \bb{G}_m(\dC)$ is called the \textit{geometric crystal with highest weight $t$}, and is identified with the open projected Richardson variety $G\overset{\circ}{/}P\subset G/P$. Moreover, $(G\overset{\circ}{/}P, f|_{X_1})$ is the mirror Landau--Ginzburg model of $\check{G}/\check{P}$ as mentioned above.

    The \textit{character $\scr{D}$-module} of the geometric crystal $X_{(G,P)}$ is defined as
    \begin{equation}\label{eq:char-D-mod}
        \mathrm{Cr}_{G,P}:=\pi_+\cE^f.
    \end{equation}
    According to \cite[Thm.~1.8]{Lam-Templier} and \cref{thm:zhu}, the character $\scr{D}$-module  $\mathrm{Cr}_{G,P}$ is a $\scr{D}$-module concentrated in degree $0$ and is smooth on $\bb{G}_m$ such that
\begin{equation}\label{eq:identification}
    \mathrm{Cr}_{G, P}\simeq \nabla_{\check{G}}(V_{\lambda_P}),
\end{equation}
where $\lambda_P$ is the minuscule weight of the dual group $\check{G}$ of $G$ determined by the parabolic subgroup $P$, and $V_{\lambda_P}$ is the representation of the highest weight ${\lambda_P}$. In particular, the fiber of $\nabla_{\check{G}}(V_{\lambda_P})$ at $1$ is the twisted de Rham cohomology $\bb{H}^n(\Omega_Y^\bullet,\rr{d}+\rr{d}w)$ associated with the Landau--Ginzburg model $(Y,w)$.

By \cref{thm:Mochizuki}, we can view the rescalable integrable mixed twistor $\scr{D}$-module $\nabla_{\check{G}}^{\rr{H}}(V)$ in \eqref{eq:EMHM-nabla} as an exponential mixed Hodge module. By \eqref{eq:identification}, its de Rham fiber is identified with $\rr{Cr}_{G,P}$. The fiber (or the pull-back) of $\nabla_{\check{G}}^{\rr{H}}(V)$ at the smooth point $1$ is an exponential mixed Hodge module (or structure) $\nabla_{\check{G}}^{\rr{H}}(V)_1$ of weight $0$, whose  de Rham fiber is $(\rr{Cr}_{G,P})_1\simeq \bb{H}^n(\Omega_Y^\bullet,\rr{d}+\rr{d}w)(n/2)$, see also \cite[Def.~A.18]{FSY22}.

By \cite[Prop.~11.22]{Mochizuki}, the irregular Hodge filtration $F_{\rr{irr}}$ of $\nabla_{\check{G}}^{\rr{H}}(V)$ induces that of $\nabla_{\check{G}}^{\rr{H}}(V)_1$. Moreover, the irregular Hodge filtration on the de Rham fiber of $\nabla_{\check{G}}^{\rr{H}}(V)_1$ coincides with the Yu filtration on the twisted de Rham cohomology \cite[Prop.~1.7.4]{Esnault-Sabbah-Yu} with a shift by $n/2$. Therefore, we have
\[
    h_{\rr{irr}}^{p-\frac{n}{2}}=\dim \rr{gr}^{p}_{F_{\rr{Yu}}} \bb{H}^{n}(\Omega_Y^\bullet,\rr{d}+\rr{d}w)=f^{p,n-p}.
\]
Lastly, since $\bb{H}^{p+q} (\Omega_Y^\bullet,\rr{d}+\rr{d}w)=0$ when $p+q\neq n$, we have $f^{p,q}=0$ when $q\neq n-p$.
\end{proof}

\begin{proof}[Proof of \upshape\cref{thm:Kontsevich}]

It is well known (see e.g.~\cite{Chevalley} or \cite[Chap.~3]{ChrissGinzburg}) that the cohomology of a homogeneous space $\check{G}/\check{P}$ is of Hodge--Tate type, i.e.~$h^{p,n-q}(X) = 0$ if $q\neq n-p$ and $h^{p,p}(X)=b_{2p}(X)$. Hence, we need to identify the Betti numbers of $X$ with our
$h^{p-\frac{n}{2}}_{irr}$. Since the latter are given by $\rho$,
it is sufficient to show that $\rr{H}^{2p}(X)$ can be identified with the eigenspace of $2\rho$ acting on $V$, which under our assumptions is a minuscule representation of $\check{G}$. This is exactly what is stated in \cite[Prop.~4.12]{Lam-Templier} (especially the last sentence of its proof). Lastly, by \cref{lem:Yu-vs-irr}, we also have the identity $h^{p,n-q}(X)=f^{p,q}(Y,w)$.
\end{proof}

\section{Examples of some small representations}\label{sec:exe}
In this subsection, we give examples of representations $V$ such that $N(V)$ have Jordan blocks of different sizes. Following our \cref{prop:Hodge-numbers}, this yields concrete results for the irregular Hodge numbers of Frenkel--Gross connections in these cases.

\subsection{\texorpdfstring{$A_{n}$}{}}\label{subsec:Examples}
We take $V=\dC^{n+1}$ as the standard representation of $\check{G}=\rr{SL}_{n+1}$. In this case the connection $\nabla_{\check{G}}(V)$ is the connection corresponding to the Bessel differential equation given in \cref{exe:An}, and we observe that $N(V)$ is conjugate to a Jordan block of size $n+1$ with eigenvalues $0$. So the local monodromy of $\nabla_{\check{G}}(V)$ at $0$ consists of $1$ single Jordan block, and the Hodge numbers are
    \[
        h^{\alpha}=
        \begin{cases}
            1, & \alpha \in \{-\frac{n+1}{2}+i\mid 0\leq i\leq n+1\},\\
            0, & \text{else},
        \end{cases}
    \]
for $\alpha\in \frac{1}{2}\dZ$.
Our result does agree with the irregular Hodge numbers of the Kloosterman connection; see for example, \cite{CDRS19,sabbah-hypergeometric,QX23}.

For the standard representation, the corresponding flag variety is the projective space $\bb{P}^n$. By \cref{thm:Kontsevich}, we deduce that well-known formula of the  Hodge numbers of $\bb{P}^n$
\[
    \{h^{p,p}\mid 0\leq p\leq n\}=\{1,\dots,1\}.
\]

\subsection{\texorpdfstring{$B_n$}{}}
Let $V=\dC^{2n+1}$ be the standard representation of $\rr{SO}_{2n+1}$. By  \cite[\S 6.3 Equation (5)]{FrenkelGross}, $N(V)$ is conjugated to a matrix with one Jordan block of size $2n+1$. So the local monodromy of $\nabla_{\check{G}}(V)$ at $0$ consists of a single Jordan block of sizes $2n+1$ and the Hodge numbers are
    \[
        h^{\alpha}=
        \begin{cases}
            1, & |\alpha|\leq n-1,\\
            0, & \text{else},
        \end{cases}
    \]
for $\alpha\in\dZ$.

This also agrees with known examples. From our result for $\nabla_{\check{G}}(V)$, we deduce the irregular Hodge numbers of $\nabla_{\rr{SO}_{2n+1}}(V)$, which coincide with formulas given in \cite{CDS21,sabbah-hypergeometric,QX23} up to a shift.

\subsection{\texorpdfstring{$E_6,F_4$}{}}
Let $V$ be a minuscule representation of $E_6$, and $\widetilde V$ the one of $F_4$. Recall that we have shown that $\nabla_{E_6}(V)=\nabla_{F_4}(\widetilde V)\oplus \cal{O}$ in \eqref{eq:decomposition-f4-e6}.

To analyze the Jordan blocks of $N(V)$, we turn to the graph \cite[(2.25)]{SW23}. Each vertex represents a basis vector $v_i,v'_i$, or $v_i''$ of $V$. An edge with number $k$ between a vertex $u:=v_i,v'_i$, or $v_i''$ with anther vertex $w:=v_{i+1},v'_{i+1}$, or $v_{i+1}''$ means that $X_{\alpha_k}u$ is a non-zero multiple of $w$. Using the representation theory of $\mathfrak{sl}_2$, we deduce that the Jordan blocks of $N(V)$ are of sizes $19$, $7$, and $1$ respectively. So the irregular Hodge numbers of $\nabla_{F_4}(\widetilde V)$ and $\nabla_{E_6}(V)$ are
\[
    h^\alpha=\begin{cases}
        1, &  4< |\alpha|\leq 8,\\
        2, &  |\alpha|\leq 4,\\
        0, &\text{else}
    \end{cases}
    \quad \text{and} \quad
    h^\alpha=\begin{cases}
        1, &  4< |\alpha|\leq 8,\\
        2, &  0<|\alpha|\leq 4,\\
        3, &  \alpha=0,\\
        0, &\text{else},
    \end{cases}
\]
respectively, for $\alpha\in\dZ$.

For the group $E_6$, the flag variety corresponds to the minuscule representation $V$ is the Cayley plane. By \cref{thm:Kontsevich}, we deduce   Hodge numbers of the Cayley plane are
\[
    \{h^{p,p}\mid 0\leq p\leq 16\}=\{1,1,1,1,2,2,2,2,3,2,2,2,2,1,1,1,1\}.
\]

\subsection{\texorpdfstring{$E_7$}{}}
Let $V$ be a minuscule representation of $E_7$. By the graph \cite[(2.27)]{SW23}, we deduce similarly as above from the representation theory of $\mathfrak{sl}_2$ that the Jordan blocks of $N(V)$ are of sizes $28$, $18$, and $10$, respectively. So the irregular Hodge numbers of $\nabla_{E_7}(V)$ are
\[
    h^{p}=\begin{cases}
        1, &  9< |p|\leq 14,\\
        2, &  5< |p|\leq 9,\\
        3, &  |p|\leq 5,\\
        0, &\text{else},
    \end{cases}
\]
for $p\in \frac{1}{2}+\dZ$.

For the group $E_7$, the flag variety corresponds to the minuscule representation $V$ is the Freudenthal variety. By \cref{thm:Kontsevich}, we deduce that the Freudenthal variety has Hodge numbers \[
    \{h^{p,p}\mid 0\leq p\leq 27\}=\{1,1,1,1,1,2,2,2,2,3,3,3,3,3,3,3,3,3,3,2,2,2,2,1,1,1,1,1\}.
\]

\subsection{\texorpdfstring{$E_8$}{}}

Let $V$ be the adjoint representation $V=\mathfrak{e}_8$ of $E_8$ (of dimension $248$). Notice that $N\in\mathfrak{e}_8$ is principally nilpotent, which can be enhanced into a principal $\mathfrak{sl}_2$. By \cite[\S6.5 \& Cor.~8.7]{Konstantprincipalthreedimensional59}, we can decompose $V=\mathfrak{e}_8$ into $\bigoplus_{i} V_{k_i}$, where $\dim V_{k_i}$ are $(2k_i+1)$-dimensional irreducible $\mathfrak{sl}_2$-representations, and the numbers $n_1< \ldots < n_8$ are the \textit{exponents} of $\mathfrak{e}_8$.

By \cite[Place VII]{BourbakiBourbakiLiegroups02}, the exponents of $\mathfrak{e}_8$ are $N_1=1$, $N_2=7$, $N_3=11$, $N_4=13$, $N_5=17$, $N_6=19$, $N_7=23$, and $N_8=29$ respectively. So the dimensions of irreducible representations (as well as the Jordan blocks of $N$) are of sizes $2N_i+1$ for $i=1,\dots,8$ respectively. So the irregular Hodge numbers of $\nabla_{E_8}(V)$ are
\[
    h^{p}=\begin{cases}
        8, & |p|\leq 1\\
        9-i, & N_{i-1}<|p|\leq N_i \text{ for } 2\leq i\leq 8\\
        0, &\text{else}
    \end{cases}
\]
for $p\in\dZ$.
\section{The \texorpdfstring{$\cR$}{R}-module associated to the Frenkel-Gross connection} \label{sec:RMod}

We have established in the previous sections using some geometric arguments that the Frenkel--Gross connection $\nabla_{\check{G}}$ can be upgraded to a tensor functor
$\nabla^{\mathrm{H}}_{\check{G}} \colon \mathrm{Rep}({\check G})\to \mathrm{MTM}^{\rr{int}}_{\rr{resc}}(\G_{m,t})$
(see the discussion around \eqref{eq:EMHM-tensor}).
Hence, for any $V\in\mathrm{Rep}(\check{G})$, the object $\nabla^{\mathrm{H}}_{\check{G}}(V)$ can in particular be considered as an integrable mixed twistor module on $\G_{m,t}$. Therefore, letting $\cR:=\cR^{\rr{int}}_{\G_{m,t}}$ be the sheaf of rings with global sections equal to $\dC[z,t^\pm]\langle z^2\partial_z,z\partial_t\rangle$, we have a $\cR$-triple $(\cM,\cM',C)$, where $\cM,\cM'$ are coherent $\cR$-modules. Our conjecture concerns an explicit expression for $\cM$.
It seems difficult to establish this conjecture directly, however, if one could do so, the main result of this paper would follow almost immediately.

For an almost simple group ${\check G}$, we associate a tensor functor
\[
\nabla_{\check{G}}^\cR \colon \mathrm{Rep}({\check G})\to \mathrm{Mod}(\cR)
\]
as follows. For each finite-dimensional complex representation $V$ of ${\check G}$, let $\cE(V):=\cO_{\dA^1_z\times \G_{m,t}}\otimes_\dC V$ be the trivial bundle on $\dA^1_z\times\mbg_{m,t}$, the associated connection on $\cE(V)$ is given by
\begin{equation}\label{eq:integrable_connection}
    \nabla_{\check{G}}^\cR(V)=\dif+\bigl(N(V)+t E(V)\bigr)\frac{\dif t}{tz} - \bigl(N(V)+tE(V)\bigr)h(\check{G})\frac{\dif z}{z^2}+\rho(V)\frac{\mathrm{d}z}{z},
\end{equation}
where $h({\check G})$ is the Coxeter number of ${\check G}$, $N(V)$ and $E(V)$ are the endomorphisms of $V$ induced by the action of $N, E\in\check{\g}$ and $\rho(V)$ is the semisimple (i.e.~diagonalizable) matrix induced by the action of $\rho=\frac12\sum_{\al^\vee\in(\Phi^\vee)^+}\al^\vee$, half sum of positive coroots of $\check G$. \\

Notice that we have
\[
\nabla_{\check{G}}^\cR(V):\cE(V)\to\cE(V)\otimes \frac{1}{z}\Omega^1_{\dA^1_z\times\G_{m,t}}\left(\log (\{0\}\times\G_{m,t})\right).
\]

\begin{proposition}\label{prop:Integrability}
    The pair $(\cE(V),\nabla_{\check{G}}^\cR(V))$ is integrable, i.e.~we have $(\nabla_{\check{G}}^\cR(V))^2=0$. In particular, the localization $(\cE(V)(*(\{0\}\times\G_{m,t})),\nabla_{\check{G}}^\cR(V))$ is a flat meromorphic connection, and yields a local system on $\mbg_{m,z}\times\mbg_{m,t}$.
\end{proposition}

\begin{proof}
First, recall that a connection $\nabla$ of the form $\nabla = A\dif t+B\dif z$ (for $A,B\in\textup{Mat}(\dim(V)\times\dim(V),\cO_{\G_{m,z}\times\G_{m,t}})$) is integrable if and only if $A$ and $B$ satisfy $\p_zA-\p_tB=[A,B]$. The connection in \eqref{eq:integrable_connection} is of this form with $A=\frac1{tz}(N+tE)$ and $B=-\frac h{z^2}(N+tE) + \frac1z\rho$. The left side of the condition is easy enough to evaluate: $\p_zA-\p_tB = \frac1{tz^2}\bigl(-N+t(h-1)E\bigr)$.

For the right side, we work inside the Lie algebra: $[A, B] = \frac1{tz^2}\bigl[N+tE,\rho \bigr]$. Recall that we have shown in \eqref{eq:[N,rho]=-N} that
    $[N,\rho]=-N. $
On the other hand, it is well-known that $\langle\theta,\rho\rangle = h-1$, so we find that $[E,\rho]=(h-1)E$.

Hence, we conclude that $[A,B]=\frac1{tz^2}\bigl(-N+t(h-1)E\bigr) = \p_zA-\p_tB$, i.e.~that $\nabla$ is integrable.
\end{proof}

We finish this paper by the following conjecture about an explicit expression of the irregular Hodge module structure defined on the Frenkel-Gross connection.
\begin{conjecture}\label{conj:IrrHodgeFilt}
For any $V\in \mathrm{Rep}({\check G})$, if $(\cM,\cM',C)$ is the $\cR$-triple of the rescalable integrable mixed twistor module $\nabla_{\check{G}}^{\rr{H}}(V)$, then we have an isomorphism $\cM\cong \nabla_{\check{G}}^\cR(V)$ of integrable $\cR$-modules.
Consequently, the irregular Hodge filtration of $\nabla_{\check{G}}^{\rr{H}}(V)$ is determined by $\rho$ in the sense of \cref{def:FiltGivenCocarac}.
\end{conjecture}
One can verify this conjecture in the examples from \cref{subsec:Examples} coming from classical groups. As mentioned above, when assuming the conjecture, our main result \cref{thm:main} would be a direct consequence using a similar strategy to \cite[Prop.~4.6]{CDS21} or to \cite[Thm.~5.9]{CDRS19} (compare the expression \eqref{eq:integrable_connection} to the expression of the connection $\nabla$ in \cite[Lem.~4.3]{CDS21}), where an adapted basis for the irregular Hodge filtration was constructed. It is, however, unclear how to show the above conjecture when the irregular Hodge structure $\nabla_{\check{G}}^{\rr{H}}(V)$ is defined using the geometric construction from \cref{sec:HNYconstruction}.

\begin{remark} Recently, Yi demonstrated in \cite{yi2025} that the $\check{G}$-Airy connections coincide with those constructed via the geometric Langlands program by Jakob-Kamgarpour-Yi in \cite{JakobKamgarpourYi}. This suggests that a similar conjectural framework as in \cref{conj:IrrHodgeFilt} may be applicable to $\check{G}$-Airy connections. In particular, one can attempt to construct, as above, $\cR$-modules associated with these connections and conjecture that they form mixed twistor $\scr{D}$-modules. With explicit descriptions of such $\cR$-modules, one may then extract the corresponding irregular Hodge filtrations. It is, however, currently even unclear how to construct an $\cR$-module from the $\check{G}$-Airy connections that satisfies the integrability property shown for the Frenkel--Gross case in \cref{prop:Integrability}.

Moreover, the methods from the previous sections used to compute the irregular Hodge numbers for Frenkel–Gross connections do not extend to the case of $\check{G}$-Airy connections, which are regular at the origin. Consequently, new techniques need to be developed to investigate the irregular Hodge theory in this setting.
\end{remark}
\bibliographystyle{amsalpha}
\providecommand{\bysame}{\leavevmode\hbox to3em{\hrulefill}\thinspace}
\providecommand{\MR}{\relax\ifhmode\unskip\space\fi MR }
\providecommand{\MRhref}[2]{  \href{http://www.ams.org/mathscinet-getitem?mr=#1}{#2}
}
\providecommand{\href}[2]{#2}

\vspace*{1cm}

\nd
Yichen Qin\\
School of Mathematical Sciences\\
Fudan University\\
Handan Road 220, Shanghai 200437 \\
China\\
yichenqin@fudan.edu.cn\\

\nd
Christian Sevenheck\\
Fakult\"at f\"ur Mathematik\\
Technische Universit\"at Chemnitz\\
09107 Chemnitz\\
Germany\\
christian.sevenheck@mathematik.tu-chemnitz.de\\

\nd
Peter Spacek\\
Max-Planck-Institute for Mathematics in the Sciences\\
04103 Leipzig\\
Germany\\
peter.spacek@mis.mpg.de\\

\end{document}